\documentclass[sn-mathphys-ay]{sn-jnl}
\bibpunct{(}{)}{;}{a}{ }{,}
\usepackage{amsmath,amssymb}
\usepackage{amsthm}

\usepackage{graphicx,color}

\usepackage{color}

\usepackage{bm}

\makeatletter

\newcommand{\Proofname}{Proof}
\makeatother

\newtheorem{theorem}{Theorem}[section]
\newtheorem{lemma}[theorem]{Lemma}

\newtheorem{proposition}[theorem]{Proposition}
\newtheorem{Definition}[theorem]{Definition}
\newtheorem{Rem}[theorem]{Remark}
\newenvironment{definition}{\begin{Definition}\rm}{\end{Definition}}
\newenvironment{remark}{\begin{Rem}\rm}{\end{Rem}}

\newcommand\be{\begin{equation}}
\newcommand\ee{\end{equation}}
\newcommand\bea{\begin{eqnarray}}
\newcommand\eea{\end{eqnarray}}
\newcommand\beaa{\begin{eqnarray*}}
\newcommand\eeaa{\end{eqnarray*}}
\newcommand\bay{\begin{array}}
\newcommand\eay{\end{array}}
\newcommand\ba{\begin{align}}
\newcommand\ea{\end{align}}
\newcommand\beba{\begin{equation}\left\{\begin{array}{rcl}}
\newcommand\eeba{\end{array}\right.\end{equation}}
\newcommand\bebaa{\begin{equation*}\left\{\begin{array}{rcl}}
\newcommand\eebaa{\end{array}\right.\end{equation*}}
\newcommand\beca{\begin{equation}\left\{\begin{array}{rcll}}
\newcommand\eeca{\end{array}\right.\end{equation}}
\newcommand\becaa{\begin{equation*}\left\{\begin{array}{rcll}}
\newcommand\eecaa{\end{array}\right.\end{equation*}}

\newcommand{\eps}{\varepsilon}

\newcommand{\R}{\mathbb{R}}
\newcommand{\Z}{\mathbb{Z}}
\newcommand{\N}{\mathbb{N}}

\newcommand{\red}{}

\def\vector#1{\mbox{\boldmath$#1$}}

\def \la{\left\langle}
\def \ra{\right\rangle}
\def \ld{\left\|}
\def \rd{\right\|}
\def \lr{\left(}
\def \rr{\right)}

\begin{document}

\title[Haptotaxis and Chemotaxis]{Relationship between haptotaxis and chemotaxis in cell dynamics 
}

\author[1]{\fnm{Hiroshi} \sur{Ishii}}\email{hiroshi.ishii@es.hokudai.ac.jp}

\author*[2]{\fnm{Hideki} \sur{Murakawa}}\email{murakawa@math.ryukoku.ac.jp}

\author[3]{\fnm{Yoshitaro} \sur{Tanaka}}\email{yoshitaro.tanaka@gmail.com}


\affil[1]{\orgdiv{Research Center of Mathematics for Social Creativity, Research Institute for Electronic Science}, \orgname{Hokkaido University}, \orgaddress{\street{Kita 12, Nishi 7, Kita-ku}, \city{Sapporo}, \postcode{060-0812}, \state{Hokkaido}, \country{Japan}}}

\affil*[2]{\orgdiv{Faculty of Advanced Science and Technology}, \orgname{Ryukoku University}, \orgaddress{\street{1-5 Yokotani, Seta Oe-cho}, \city{Otsu}, \postcode{520-2194}, \state{Shiga}, \country{Japan}}}

\affil[3]{\orgdiv{Department of Complex and Intelligent Systems, School of Systems Information Science}, \orgname{Future University Hakodate}, \orgaddress{\street{116-2 Kamedanakano-cho}, \city{Hakodate}, \postcode{041-8655}, \state{Hokkaido}, \country{Japan}}}

\abstract{
    { Cell sorting mediated by direct cell--cell contact or cellular protrusions can be described by nonlocal cell--cell adhesion models of haptotaxis type, whereas communication through diffusible chemical signals is described by Keller--Segel type chemotaxis systems. We investigate the mathematical relationship between these two descriptions. We first prove that, in the fast signal diffusion limit, subsequences of weak solutions of a Keller--Segel type parabolic--parabolic chemotaxis system with possibly degenerate nonlinear diffusion and density saturation converge to weak solutions of the corresponding parabolic--elliptic system. Eliminating the elliptic chemical fields rewrites the population equation as a nonlocal adhesion model whose interaction kernel is a finite linear combination of Green functions. We then prove, in arbitrary spatial dimensions, that the periodization of the gradient of a prescribed radially symmetric interaction potential can be approximated by finite linear combinations of gradients of Green functions. Combining these results, we prove that suitable weak solutions of the nonlocal adhesion model and a parabolic--parabolic chemotaxis system can be chosen with an arbitrarily small difference over any fixed finite time interval. Numerical simulations compare the nonlocal, parabolic--elliptic, and parabolic--parabolic models and examine how the number of kernel terms and the relaxation time affect their differences.}
}

%

\maketitle

\section{Introduction}

{
Cell sorting is the spatial rearrangement of initially intermixed cells into organized populations or tissue structures. It is driven by several mechanisms, including cell--cell adhesion and communication through diffusible chemical signals, which may act simultaneously during development.

Direct cell--cell contact and protrusion-mediated sensing motivate nonlocal cell--cell adhesion models in which migration depends on surrounding cell densities over a finite sensing range. Such models provide haptotaxis-type descriptions of contact-mediated migration, whereas chemotaxis describes migration in response to gradients of diffusible chemicals. In this paper, we compare a nonlocal cell--cell adhesion model, hereafter called the nonlocal haptotaxis model, with Keller--Segel type chemotaxis systems.
}

{
\citet{armstrong2006} introduced a nonlocal continuum model of cell--cell adhesion in which the adhesive velocity is given by an integral of the surrounding cell density. \citet{murakawaTogashi2015} subsequently reformulated the basic cell movement in terms of population pressure. Building on these models, \citet{carrillo2019population} incorporated population pressure and density saturation, leading to the following model:
}
\begin{equation}
\label{eq:adh}
	\frac{\partial  u}{ \partial t} = \nabla \cdot (u\nabla \chi (u))  - \nabla \cdot (g(u){\vector{V}}(u)),
\end{equation}
where $u=u({x},t)$ is the population density of cells at position $x$ and time $t$. The first term of the {right-hand side} of \eqref{eq:adh} describes the movement of cells due to pressure, where $\chi(u)$ denotes the pressure. 
The second term of the {right-hand side} of \eqref{eq:adh} represents the movement of cells due to adhesion. Here, the velocity is defined by 
\begin{align*}
{\vector{V}}(u)({x}) 
=  \int_0^1\!\!\!  \int_{S^{N-1}} u({x}+r{\eta})\, \omega(r)r^{N-1} {\eta} \, d{\eta} dr , 
\end{align*}
where $S^{N-1}$ is { the $(N-1)$-dimensional unit sphere}. 
The function $\omega$ describes how the force is dependent on the distance from ${x}$. 
This controls cell adhesion and repulsion. Accordingly, {$\vector{V}(u)$} implies that each cell counts its surrounding within the sensing radius, that is rescaled to be one, to determine the direction of movement. 
{ For a pressure normalized by $\chi(1)=1$, a typical example of $g$ is $g(u)=u(1-\chi(u))$ for $0\leq u\leq 1$, with $g(u)=0$ outside $[0,1]$. This indicates that the magnitude of the total force of adhesion and/or repulsion decreases as the pressure increases with the local density and vanishes at the maximal density $u=1$. This term restricts the solution to lie between $0$ and $1$. This effect is called the density saturation effect or volume filling effect~\citep{HillenPainter2001,painterHillen2002}.} 

{ This model has been extended to a multi-component version that reproduces cell sorting driven by cell--cell adhesion~\citep{carrillo2019population}.} Moreover, it has been applied to solving real-world problems in life sciences~\citep{carrillo2024ara,matsunaga2017reelin,trush}.

Shifting the discussion to chemotaxis, the movement of cell populations due to chemotaxis has been studied using the Keller--Segel model{~\citep{HillenPainter2009}}. Here, we consider the following problem, where the fundamental movement of cells is described by the aforementioned pressure-dependent diffusion, and chemotaxis is represented by a Keller--Segel type system with a density saturation effect.

\begin{equation}
\label{eq:KS}
\left\{
\begin{aligned}
	\frac{ \partial u}{ \partial t } & =   \nabla \cdot \left(u\nabla \chi (u) \right)  - \nabla \cdot (g(u)\nabla (av)),\\
	\xi \frac{ \partial v }{ \partial t } &= d \Delta v - v + u, 
\end{aligned}
\right.  
\end{equation}
where $u=u({x},t)$ and $v=v({x},t)$ are the population density of cells and the concentration of the chemical substance, respectively. 
The constant $d$ {denotes} the diffusion coefficient of the chemical substance and { $\xi$ is the relaxation time of the chemical dynamics; a smaller $\xi$ corresponds to faster equilibration of the chemical field.}  
The parameter $a\in \R$ represents the strength of chemotaxis, with 
$a>0$ indicating that the chemical acts as an attractant, and $a<0$ indicating that the chemical acts as a repellent. 
The second equation in \eqref{eq:KS} implies that the chemical diffuses, degrades, and is produced by the cells. 

Both the haptotaxis model~\eqref{eq:adh} and the chemotaxis model~\eqref{eq:KS} represent cell dispersion and aggregation/repulsion, but there are also significant differences. { When $\chi(u)=u^\gamma$ ($\gamma >0$)}, the first term on the right-hand side of both \eqref{eq:adh} and \eqref{eq:KS} is a porous medium type nonlinear diffusion, which has the property of finite propagation of the solution. This means that if the support of the initial value is compact, the support of the solution also remains compact. Furthermore, the aggregation term in \eqref{eq:adh} only involves interactions within the sensing radius. As a result, cell populations following \eqref{eq:adh} can form multiple colonies. { On the other hand, in system~\eqref{eq:KS} with $a>0$, the chemical field has global support, so long-range attraction may promote coarsening toward fewer aggregates. Despite these differences, the two models can exhibit similar behavior over finite time intervals. The present results concern the relationship between the two models over finite time intervals and do not address long-time pattern selection.}

Investigating these differences is important for understanding cell aggregation/repulsion, {cell--cell} adhesion, and cell sorting. { In this work, we investigate the relationship between the solution behavior of \eqref{eq:adh} with a prescribed interaction kernel and that of a Keller--Segel type chemotaxis system.}


We consider problems in an $N$-dimensional bounded domain $\Omega=[-L,L)^N$ ($L>0$) with periodic boundary condition. Therefore, $\Omega$ coincides with $\R^N / (2L\Z)^N$. 
Defining $\beta(u)=\int_0^u s\chi' (s)\, ds$ and { extending $\omega$ by zero outside $[0,1]$, let $\nabla K({x})=-\omega(|{x}|) \frac{{x}}{|{x}|}$ denote the whole-space interaction field, and let $\nabla W$ denote its periodization. Then} \eqref{eq:adh} can be rewritten as follows:
\begin{equation}
\label{eq:AD}
	\frac{\partial  u}{ \partial t} = \Delta \beta(u)  - \nabla \cdot \left(g(u) \nabla W*u \right) \qquad \mathrm{in }\ Q_T:=\Omega \times (0,T]. 
\end{equation}
{ In particular, $\chi(u)=u^\gamma$ yields $\beta(u)=\frac{\gamma}{\gamma+1}u^{\gamma+1}$.}
{ Here, $T>0$ is fixed, }
$*$ denotes the convolution of two periodic functions in the space variable, namely, 
\begin{equation*}
	(W*u) ({x},t) = \int_\Omega W({x}-{y}) u({y},t) d{y}.
\end{equation*}
As usual, we impose the initial condition $u(0) = u_0$. 

{ When a sequence of functions $\{W_m\}_{m\in \N}$ approximates $W$ in a suitable sense, it is natural to compare solutions of \eqref{eq:AD} with those of the following equation:}
\begin{equation}
\label{eq:ADtilde}
	\frac{\partial  u}{ \partial t} = \Delta \beta(u)  - \nabla \cdot \left(g(u) \nabla {W}_m*u \right). 
\end{equation}

If $w$ is the fundamental solution of
\begin{equation}
\label{eq:funda}
-d\Delta w +w = \delta, 
\end{equation}
the solution of \eqref{eq:ADtilde} with ${W}_m=aw$ coincides with the solution of \eqref{eq:KS} with $\xi=0$. 
Extending this, we consider the case where ${W}_m$ can be expressed as a linear combination of the fundamental solutions of \eqref{eq:funda} associated with different $d$. In this case, \eqref{eq:ADtilde} corresponds to the following problem: 
\begin{equation} 
\label{eq:KSPE} 
\left\{
\begin{aligned}
	\frac{\partial  u}{ \partial t} &= \Delta \beta(u)  - \nabla \cdot \left(g(u) \nabla \sum_{j=1}^M a_j v_j \right), \\
	0 &=   d_j \Delta v_{j}   - v_j + u  \quad (j=1,2,\dots, M).
\end{aligned}
\right.  
\end{equation}
In this model, cells produce multiple diffusive and degradable signals, which in turn drive the cells themselves. 
The constant $d_j >0$ is the diffusion coefficient, and $a_j \in \R$ {denotes} the sensitivity of the cells to the $j$-th chemical. 
{ For the parabolic--elliptic system~\eqref{eq:KSPE}, only the initial condition $u(0)=u_0$ is imposed; for each $t$, $v_j$ is determined from $u(t)$ by the elliptic equation. }

{ If the interaction field $\nabla W$ can be approximated by a linear combination of gradients of fundamental solutions, it is natural to compare weak solutions of \eqref{eq:KSPE} and \eqref{eq:AD}. In fact, we show that suitable weak solutions of the two problems can be chosen so that their $L^2(Q_T)$ difference is arbitrarily small over a finite time interval.}

The derivation of this relationship between \eqref{eq:AD} and \eqref{eq:KSPE} is motivated by the work of \citet{ninomiya2017rda,ninomiya2018rda}. They studied a method to approximate solutions of reaction-diffusion equations with nonlocal reactions using solutions of a FitzHugh--Nagumo type system composed only of local terms. In their study, they have shown that any potential $W$ can be approximated in $L^2$ space by a linear combination of fundamental solutions of certain equations in one spatial dimension. 
{
They addressed nonlocal reactions in one spatial dimension.
%
More recently, \citet{IshiiTanaka2026} extended this approach to arbitrary spatial dimensions. 
They showed that any radial kernel in
$L^1(\R^N)$ can be approximated in $L^1(\R^N)$ by a finite linear combination of fundamental solutions. 
Furthermore, under regularity and decay assumptions depending on the spatial dimension, they obtained constructive approximations in $W^{1,1}(\R^N)$ and derived quantitative error estimates for $N=1,2,3$. 
Their study concerns radial kernels on $\R^N$, whereas we focus on nonlocal advection involving the gradient of a periodic interaction potential.
}

{ System~\eqref{eq:KSPE} treats the chemical fields as quasi-static. To incorporate a finite response time of the chemical fields, we consider the following parabolic--parabolic system:}
\begin{equation}
\label{eq:KSM} 
\left\{
\begin{aligned}
	\frac{\partial  u}{ \partial t} &= \Delta \beta(u)  - \nabla \cdot \left(g(u) \nabla \sum_{j=1}^M a_j v_j \right), \\
	\xi \frac{ \partial v_j }{ \partial t } &=   d_j \Delta v_{j}   - v_j + u  \quad (j=1,2,\dots, M).
\end{aligned}
\right.  
\end{equation}
Here $0<\xi \le 1$ is a relaxation parameter.  
{ For the parabolic--parabolic system~\eqref{eq:KSM}, we impose the initial condition
\begin{equation}
\label{init_uv}
(u(0), v_1(0), \dots , v_M(0)) = (u_0, v_{01},\dots , v_{0M}).     
\end{equation}
}

\citet{murakawa_tanaka2024} examined the relationship between \eqref{eq:AD} and \eqref{eq:KSM} in the case where the diffusion is linear ($\beta(u)=u$) and there is no density saturation effect ($g(u) = u$) in one spatial dimension.
{ They} demonstrated that the solutions of these two equations are close if $M\in \N$, $\{ a_j \}_{j=1}^M$, $\{ d_j \}_{j=1}^M$ and $\xi >0$ are chosen appropriately. Furthermore, using this Keller--Segel type approximation, { they} showed that the destabilization of solutions near equilibrium points of the nonlocal Fokker-Planck equation is very similar to diffusion-driven instability. 
In this work, we analyze the cases with nonlinear diffusion and the density saturation effect in multiple spatial dimensions. 

{
Nonlocal-to-local limits have recently been studied for aggregation--diffusion equations with nonlocal cell--cell interactions, where local equations are derived by shrinking the interaction range. \citet{falcoBakerCarrillo2024} derived a local continuum model from a nonlocal model of cell--cell adhesion, and \citet{falcoBakerCarrillo2025} discussed a nonlocal-to-local approach to aggregation--diffusion equations. Moreover, \citet{carrilloElbarSkrzeczkowski2025} and \citet{elbarSkrzeczkowski2023} analyzed nonlocal-to-local limits for degenerate Cahn--Hilliard equations and systems.

In contrast to these nonlocal-to-local limits, we approximate the gradient $\nabla W$ of the interaction potential appearing in the nonlocal haptotaxis model, namely the haptotactic interaction kernel, by a finite linear combination of gradients of Green functions with different diffusion coefficients. This provides a mathematical framework for relating the nonlocal haptotaxis model to a Keller--Segel type chemotaxis system with finitely many chemical fields.
}

{
The contributions of this paper are twofold. First, we provide a rigorous framework connecting a nonlocal haptotaxis model with a parabolic--parabolic Keller--Segel type system containing a finite number of chemical fields. By combining the fast signal diffusion limit with a kernel approximation, we establish the convergence of the corresponding weak solutions. This extends the previous one-dimensional result for linear diffusion without density saturation~\citep{murakawa_tanaka2024} to multiple spatial dimensions with possibly degenerate nonlinear diffusion and density saturation. To the best of our knowledge, the relationship between haptotaxis and chemotaxis has not previously been established from a parabolic--parabolic system in this setting. Second, in Section~\ref{sec:approx_kernels}, we show that, in arbitrary spatial dimensions, the gradient of a given radially symmetric interaction potential can be approximated by finite linear combinations of gradients of fundamental solutions. This kernel approximation is also of independent mathematical interest.}

{ In addition, we numerically compare the solutions of the three models on a one-dimensional periodic domain and examine how the number $M$ of terms in the kernel approximation and the relaxation time $\xi$ affect the differences among the numerical solutions.}

{ The organization of this paper is as follows. Section~\ref{sec:results} introduces the assumptions and states the main results. Section~\ref{sec:existence} establishes the existence of weak solutions for the Keller--Segel type system, and Section~\ref{sec:pptope} proves the fast signal diffusion limit. Section~\ref{sec:approx_kernels} addresses the approximation of gradients of interaction potentials, while Section~\ref{sec:convWm} establishes convergence under kernel approximation and proves the main comparison result. Section~\ref{sec:numerical} numerically compares the three models and examines the effects of $M$ and $\xi$ on the differences among their solutions. The final section presents concluding remarks.}

\section{Assumptions and main results}
\label{sec:results}

In this work, we impose the following assumptions. 
\renewcommand{\theenumi}{H\arabic{enumi}}
\renewcommand{\labelenumi}{(\theenumi)}
\begin{enumerate}
\item \label{assum_gamma}
$\beta\in C^2(\R)$ is { a} strictly increasing function with $\beta(0)=0$. 
\item \label{assum_g}
{ $g\in W^{1,\infty}(\R)$ with Lipschitz constant $L_g$ satisfies $g(s)=0$ for $s\notin[0,1]$. 
}
\item \label{assum_init_u}
$0\le u_0 \le 1$ a.e. 
\item \label{assum_init_v}
{ $v_{0j}\in W^{1,\infty}(\Omega)$ ($j=1,2,\dots, M$).}
%
%
%
%
\item \label{assum_W}
$W$ is a periodic function satisfying
\begin{equation}
\label{wholespace_to_periodic}
    {
    \nabla W=\sum_{\ell\in \Z^N} \nabla K(\,\cdot-2L\ell) \quad \mbox{ in } L^1(\Omega)
    }
\end{equation}
for a radial function {$K$ satisfying $|\nabla K|\in L^{1}\lr \R^N\rr$}.
\end{enumerate}
{
For haptotaxis models on $\mathbb{R}^N$, finite-range interactions are often
considered; in the present potential formulation this corresponds to compact
support of the interaction field $\nabla K$. Under periodic boundary conditions,
we therefore periodize the interaction field, rather than the potential itself, as
in \eqref{wholespace_to_periodic}. If $\nabla K$ is compactly supported, then only
finitely many terms in \eqref{wholespace_to_periodic} are nonzero, with the number
of such terms bounded uniformly in $x\in\Omega$. Assumption~\eqref{assum_W}
covers and extends this finite-range setting by requiring only
$|\nabla K|\in L^1(\mathbb{R}^N)$. Under this assumption, the series defining
$\nabla W$ converges absolutely in $L^1(\Omega)$, and hence noncompactly
supported interaction fields are also allowed. This is sufficient because the
equations and estimates involve only $\nabla W$. 
For example, this assumption includes Gaussian interaction potentials, which are not compactly supported but have integrable gradients. It also covers the interaction fields generated by solutions of $-d\Delta v+v=u$, since the gradient of the corresponding whole-space Green function is integrable and its periodization gives the periodic interaction field. 
}

Problems \eqref{eq:AD} and \eqref{eq:KSM} are understood in the following weak sense. 
\begin{definition}
\label{weak_AD}
A function $u \in L^\infty(Q_T) \cap H^1(0,T;H^1(\Omega)^*)$ {satisfying $0\le u\le1$ a.e.\ in $Q_T$} is said to be a weak solution of { the initial-boundary value problem} \eqref{eq:AD} if it fulfills
\begin{align*}
&\beta(u) \in L^2(0,T;H^1(\Omega)), \\
&\int_0^T \la \frac{\partial u}{\partial t}, \varphi \ra 
+\int_0^T \la \nabla \beta(u),\nabla \varphi\ra
 -  \int_0^T \la g(u) \nabla W*u,\nabla \varphi\ra
 =0 \\
 & 
 \hspace*{4.5cm} \mbox{{ for all functions} } \varphi \in L^2(0,T;H^1(\Omega)), \\
 &u(0) = u_0. 
\end{align*}
\end{definition}
Here and hereafter, we denote by $\langle \cdot,\cdot\rangle $ both the inner product in $L^2(\Omega)$ and the duality pairing between $H^{1}(\Omega)^*$ and $H^1(\Omega)$. 
We use a simple notation $\|\cdot \|$ for the norms in both $L^2(\Omega)$ and $L^2(\Omega)^N$ spaces. 
\begin{definition}
\label{weak_KSM}
A set of functions $(u, v_{1},\dots , v_{M})$ is said to be a weak solution of { the initial-boundary value problem} \eqref{eq:KSM} if it satisfies \eqref{init_uv} and 
\begin{align*}
&{0\le u\le1 \quad \mbox{a.e.\ in }Q_T,}\\
u &\in L^\infty(Q_T) \cap H^1(0,T;H^1(\Omega)^*), \\ 
\beta(u) &\in L^2(0,T;H^1(\Omega)), \\
v_j &\in H^1(0,T;L^2(\Omega)) \cap L^2(0,T;H^2(\Omega)) \quad (j=1,2,\dots,M), \\
\int_0^T& \la \frac{\partial u}{\partial t}, \varphi \ra 
+\int_0^T \la \nabla \beta(u),\nabla \varphi\ra
 -  \int_0^T \la g(u) \nabla \sum_{j=1}^M a_j v_j,\nabla \varphi\ra
 =0 \\
 & 
 \hspace*{4cm} \mbox{{ for all functions} } \varphi\in L^2(0,T;H^1(\Omega)), \\
\xi \frac{ \partial v_j }{ \partial t } &=   d_j \Delta v_{j}   - v_j + u \quad \mbox{a.e.}
\end{align*}
\end{definition}
Weak solutions for other problems are defined similarly.

One of the main results is the following existence theorem.
\begin{theorem}
    \label{thm:exist:KSM}
    Under Assumptions \eqref{assum_gamma}--\eqref{assum_init_v}, a weak solution of \eqref{eq:KSM} exists. 
\end{theorem}
The existence of a weak solution of \eqref{eq:AD} can be established through the approximation process.
{ More precisely, the limits are taken sequentially: first $\xi\to0$ for each fixed kernel approximation $W_m$, and then $m\to\infty$.}
\begin{theorem}
    \label{thm:exist:AD}
    Under Assumptions \eqref{assum_gamma}, \eqref{assum_g}, \eqref{assum_init_u}, and \eqref{assum_W}, a weak solution of \eqref{eq:AD} exists. 
\end{theorem}

The weak solution of the {parabolic--parabolic} problem~\eqref{eq:KSM} is close to that of the {parabolic--elliptic} problem~\eqref{eq:KSPE} when $\xi$ is sufficiently small. 
\begin{theorem}
    \label{thm:conv:xi}
    Assume that \eqref{assum_gamma}--\eqref{assum_init_v} hold. 
    {
    For each $\xi\in(0,1]$, let $\lr u^{\xi,M},v_1^{\xi,M},\dots,v_M^{\xi,M}\rr$ be a weak solution of \eqref{eq:KSM}.
    Then, there exist a sequence $\{\xi_k\}_{k\in\N}\subset(0,1]$ with $\xi_k\to0$ and a weak solution $\lr u^{M}, v_1^{M}, \dots, v_M^{M} \rr$ of \eqref{eq:KSPE} such that
    \begin{align*}
&u^{\xi_k,M} \to u^M
& &  \mbox{strongly in } L^2(Q_T), \mbox{ a.e. in } Q_T, \\
& & & \mbox{weakly in }  H^1(0,T;H^1(\Omega)^*),  \\
&\beta \lr u^{\xi_k,M}\rr \to \beta(u^M) 
& & \mbox{strongly in } L^2(Q_T), \mbox{ a.e. in } Q_T,    \\
& & & \mbox{weakly in } L^2(0,T;H^1(\Omega)), \notag \\
&v_j^{\xi_k,M} \to v_j^M 
& & \mbox{strongly in } L^2(0,T;H^1(\Omega)), \mbox{ a.e. in } Q_T, \\
&g\lr u^{\xi_k,M}\rr \nabla \sum_{j=1}^M a_j v_j^{\xi_k,M} \\
&\quad\to g(u^M) \nabla \sum_{j=1}^M a_j v_j^M,
& & \mbox{strongly in } L^2(Q_T),\mbox{ a.e. in } Q_T
\end{align*}
    as $k$ tends to infinity.
    }
\end{theorem}

{
\begin{remark}
In the limit $\xi\to0$, the equations for $v_j$ become elliptic, and hence the initial data $v_{0j}$ are not retained as independent initial conditions in the limiting problem. The convergence of $v_j^{\xi_k,M}$ in Theorem~\ref{thm:conv:xi} is in a space-time $L^2$ space and does not assert convergence at $t=0$. Thus, an initial layer may occur when $v_{0j}$ is not compatible with the limiting elliptic equation. In Theorem~\ref{Thm_conv}, we use the compatible initial data $v_{0j}=w_j*u_0$ specified in \eqref{init_u_wu}.
\end{remark}
}

Furthermore, we obtain the following result regarding the convergence of weak solutions associated with the approximation of $W$.
\begin{theorem}
    \label{thm:conv:Wm}
    Assume that \eqref{assum_gamma}, \eqref{assum_g}, and \eqref{assum_init_u} hold. 
    Let $\{ W_m\}_{m\in \N}$ be a sequence of functions such that $\nabla W_m$ converges to $\nabla W$ in $L^{1}(\Omega)$ as $m$ tends to infinity. 
    {
    For each $m\in\N$, let $u_m$ be a weak solution of \eqref{eq:ADtilde} corresponding to $W_m$.
    Then, there exist a subsequence $\{m_k\}_{k\in\N}$ and a weak solution $u$ of \eqref{eq:AD} such that
    \begin{align*}
&u_{m_k} \to u
& &  \mbox{strongly in } L^2(Q_T), \mbox{ a.e. in } Q_T, \\
& & & \mbox{weakly in }  H^1(0,T;H^1(\Omega)^*),  \\
&\beta \lr u_{m_k}\rr \to \beta(u) 
& & \mbox{strongly in } L^2(Q_T), \mbox{ a.e. in } Q_T,    \\
& & & \mbox{weakly in } L^2(0,T;H^1(\Omega))
\end{align*}
    as $k$ tends to infinity.
    }
\end{theorem}

One of the key aspects of this work is that { the interaction field $\nabla W$ can be approximated by linear combinations of gradients of fundamental solutions.}  
We consider the Green functions $w_j$ of the following problem with the periodic boundary condition: 
\begin{equation*}
{-d_j \Delta w_j+ w_j =\delta.}
\end{equation*}
Here, $\delta$ denotes the Dirac delta function. 
With this notation, we establish the following result. 

\begin{theorem}
    \label{thm:W_approx}
    {Assume that \eqref{assum_W} holds for $W$. 
    Then, for any $\eps>0$, there exist $M\in \N$, a family of constants $\{a_j\}^M_{j=1}$ and a family of positive constants $\{d_j\}^M_{j=1}$ such that 
    \begin{equation*}
    \left\| \nabla W - \nabla \sum^{M}_{j=1} a_j w_j \right\|_{L^{1}\lr \Omega\rr} < \eps.
    \end{equation*}}
\end{theorem}
%
%

{
\begin{remark}
Theorem~\ref{thm:W_approx} is nonconstructive under the general assumption $|\nabla K|\in L^1(\R^N)$ because the proof is based on Wiener's Tauberian theorem (see Section~\ref{sec:approx_kernels}). 
In particular, Theorem~\ref{thm:W_approx} does not provide an explicit choice of $M$, $\{a_j\}_{j=1}^M$, or $\{d_j\}_{j=1}^M$, nor does it give a quantitative relation between $M$ and the prescribed error $\eps$. 
Constructive approximations available under stronger assumptions are discussed in Remark~\ref{rem:constructive_approx}.
\end{remark}
}

{For these $M$, $\{ a_j \}_{j=1}^M$ and  $\{d_j\}^M_{j=1}$}, we set $W_M=\sum_{j=1}^M a_j w_j$. 
{ Then, \eqref{eq:AD} with $W=W_M$ as the potential coincides with the parabolic--elliptic Keller--Segel type system~\eqref{eq:KSPE}, whose elliptic components satisfy $v_j(t)=w_j*u(t)$. To compare this parabolic--elliptic solution with a weak solution of the parabolic--parabolic system~\eqref{eq:KSM}, we choose the compatible initial data}
\begin{equation}
\label{init_u_wu}
{ (u^{\xi,M}(0), v_1^{\xi,M}(0), \dots , v_M^{\xi,M}(0)) = (u_0, w_1*u_0,\dots , w_M*u_0).}     
\end{equation}
Combining Theorems~\ref{thm:conv:xi}, \ref{thm:conv:Wm} and \ref{thm:W_approx}, we conclude with the following result. 
\begin{theorem}
\label{Thm_conv}
{Assume that \eqref{assum_gamma}--\eqref{assum_init_u} and \eqref{assum_W} are satisfied.}
{ Then, for any $\eps>0$ and $T>0$, there exist
$M\in \N$, a family of constants $\{ a_j \}^M_{j=1}$, a family of positive constants $\{d_j\}^M_{j=1}$, $\xi >0$, a weak solution $u$ of \eqref{eq:AD} with the initial datum $u_0$, and a weak solution $\lr u^{\xi,M},v_1^{\xi,M},\ldots,v_M^{\xi,M}\rr$ of \eqref{eq:KSM} with the initial datum \eqref{init_u_wu}, such that
\[
    \ld u - u^{\xi,M} \rd_{L^2(Q_T)}<\eps.
\]}
\end{theorem}

\section{Existence of a weak solution for the Keller--Segel type system}
\label{sec:existence}
In this section, we prove the existence of { a} weak solution of the Keller--Segel type system~\eqref{eq:KSM}. 
{
The proof in this section uses the same broad regularization and
Schauder fixed-point strategy as \citet{bendahmane2007} for the case
$M=1$, while treating an arbitrary fixed finite number $M\ge1$ of
chemical fields. However, several regularity statements, uniqueness
arguments, and uniform estimates in their proof are not justified
under the assumptions stated there. We therefore give a self-contained
proof for $M\ge1$, including the estimates and compactness arguments
needed for the subsequent degenerate limit and for
Section~\ref{sec:pptope}.

The proof proceeds in two steps. First, we introduce a uniformly parabolic approximation $\beta_\eta$, $\eta\in(0,1]$, and establish the existence of a weak solution to the resulting non-degenerate problem by the Schauder fixed-point theorem. We then derive estimates independent of $\eta$ and pass to the limit $\eta\to0$ to obtain a weak solution of the original degenerate problem.
}

\subsection{Non-degenerate approximation}

In the equation for $u$ in \eqref{eq:KSM}, the diffusion is generally degenerate.
{
To justify the existence of a solution to the non-degenerate problem
without imposing any growth condition on $\beta'$ outside $[0,1]$,
we introduce the auxiliary function
\[
\widetilde\beta(s):=
\begin{cases}
\beta'(0)s, & s<0,\\
\beta(s), & 0\le s\le1,\\
\beta(1)+\beta'(1)(s-1), & s>1.
\end{cases}
\]
Then $\widetilde\beta\in C^{1,1}(\mathbb R)$ is nondecreasing and
globally Lipschitz continuous, and
$\widetilde\beta=\beta$ on $[0,1]$.
For $\eta\in(0,1]$, we define the non-degenerate approximation by
\[
\beta_\eta(s):=\eta s+\widetilde\beta(s).
\]
For every $s\in\mathbb R$,
\[
\eta
\le \beta_\eta'(s)
\le \eta+\max_{r\in[0,1]}\beta'(r).
\]
In particular, $\beta_\eta'$ is globally bounded and Lipschitz continuous.
}
Consequently, we consider the following system:
\begin{equation}
\label{eq:KSrelax} 
\left\{
\begin{aligned}
	\frac{\partial  u}{ \partial t} &= \Delta \beta_\eta(u)  - \nabla \cdot \left(g(u) \nabla \sum_{j=1}^M a_j v_j \right), \\
	\xi \frac{ \partial v_j }{ \partial t } &=   d_j \Delta v_{j}   - v_j + u  \quad (j=1,2,\dots, M)
\end{aligned}
\right.  \qquad \mbox{in } Q_T
\end{equation}
with the initial condition \eqref{init_uv} and the periodic boundary condition. 
To prove the existence of a weak solution of \eqref{eq:KSrelax}, we deal with each equation one by one.
Set 
\[ 
\mathcal{K} 
:=
\left\{ u\in L^2(Q_T) \ : \ 0\le u \le 1\quad  \mathrm{ a.e.} \right\}, 
\] which is a closed { convex} subset of $L^2(Q_T)$. 
For a given $\bar{u}\in \mathcal{K}$, we consider the following parabolic equations: 
\begin{equation}
\label{eq:vj_indep} 
\left\{
\begin{aligned}
	\xi \frac{ \partial v_j }{ \partial t } &=   d_j \Delta v_{j}   - v_j + \bar{u} 
 & &\mbox{in } Q_T, \\
 v_j(0)&=v_{0j}& &\mbox{in } \Omega 
\end{aligned}
\right.
\end{equation}
for each $j=1,2,\dots, M$. 
{
These equations are linear and uniformly parabolic. Their unique solutions
are given by the variation-of-constants formula
\begin{align*}
v_j(t)
={}&e^{-t/\xi}e^{(d_jt/\xi)\Delta}v_{0j}\\
&+\frac{1}{\xi}\int_0^t e^{-(t-s)/\xi}
e^{(d_j(t-s)/\xi)\Delta}\bar{u}(s)\,ds .
\end{align*}
Standard $L^2$ parabolic regularity~\citep{ladyzenskaja} and the periodic
heat-semigroup estimates yield
\[
v_j\in H^1(0,T;L^2(\Omega))
\cap L^2(0,T;H^2(\Omega))
\cap L^\infty(0,T;W^{1,\infty}(\Omega)).
\]
In particular,
\[
\|\nabla v_j(t)\|_{L^\infty(\Omega)}
\le
\|\nabla v_{0j}\|_{L^\infty(\Omega)}
+C_\Omega d_j^{-1/2}\|\bar{u}\|_{L^\infty(Q_T)}
\]
for a.e. $t\in(0,T)$. Consequently,
\[
b:=\nabla\sum_{j=1}^M a_jv_j\in L^\infty(Q_T;\R^N).
\]
}
{With these $v_j$ fixed, we consider} the following quasilinear parabolic equation: 
\begin{equation}
\label{eq:u_indep} 
\left\{
\begin{aligned}
	\frac{\partial  u}{ \partial t} &= \Delta \beta_\eta(u)  - \nabla \cdot \left(g(u) \nabla \sum_{j=1}^M a_j v_j \right)& &\mbox{in } Q_T, \\
 u(0)&=u_{0} & &\mbox{in } \Omega. 
\end{aligned}
\right.  
\end{equation}
{Since $\beta_\eta'$ is globally bounded and Lipschitz continuous and satisfies $\beta_\eta'\ge\eta$, standard uniformly parabolic theory~\citep[Chapter~V, Theorem~6.7]{ladyzenskaja} applies to \eqref{eq:u_indep}. Hence, it has a weak solution $u\in L^\infty(0,T;L^2(\Omega))\cap L^2(0,T; H^1(\Omega))$.}
{
We next prove uniqueness for the fixed $\bar{u}$, and hence for the fixed
drift $b$. Let $u_1$ and $u_2$ be two weak solutions of
\eqref{eq:u_indep} with the same initial datum, and set
$U:=u_1-u_2$. Mass conservation gives
$\int_\Omega U(t)\,dx=0$. For a.e. $t\in(0,T)$, let $\phi(t)$ be the
periodic solution of
\[
-\Delta\phi=U,\qquad \int_\Omega\phi\,dx=0,
\]
and set $\|U\|_{-1,\mathrm{per}}:=\|\nabla\phi\|_{L^2(\Omega)}$.
Testing the difference of the two equations by $\phi$ gives
\begin{align*}
&\frac{1}{2}\frac{d}{dt}\|U\|_{-1,\mathrm{per}}^2
+\int_\Omega
\bigl(\beta_\eta(u_1)-\beta_\eta(u_2)\bigr)U\,dx\\
&\qquad=
\int_\Omega\bigl(g(u_1)-g(u_2)\bigr)b\cdot\nabla\phi\,dx .
\end{align*}
Since $\beta_\eta'(s)\ge\eta$ and $g$ is Lipschitz continuous,
\begin{align*}
&\frac{1}{2}\frac{d}{dt}\|U\|_{-1,\mathrm{per}}^2
+\eta\|U\|_{L^2(\Omega)}^2\\
&\qquad\le
L_g\|b\|_{L^\infty(Q_T)}
\|U\|_{L^2(\Omega)}\|U\|_{-1,\mathrm{per}}\\
&\qquad\le
\frac{\eta}{2}\|U\|_{L^2(\Omega)}^2
+\frac{L_g^2}{2\eta}\|b\|_{L^\infty(Q_T)}^2
\|U\|_{-1,\mathrm{per}}^2 .
\end{align*}
Since $U(0)=0$, Gronwall's inequality yields $U=0$. Thus the weak
solution of \eqref{eq:u_indep} is unique for each fixed $\bar{u}$.
Consequently, the map $\Theta(\bar u)=u$ used below is well-defined
and single-valued.
}
%
%
%
%
The solutions $v_j$ and $u$ of \eqref{eq:vj_indep} and \eqref{eq:u_indep} depend on $\eta$, $\xi$ and $\bar{u}$, but we omit these dependencies to simplify the notation. 
We provide uniform bounds with respect to $\xi$ and $\eta$.

\begin{lemma}
\label{Lem_linf_u}
The following boundedness result holds for the weak solution $u$ of \eqref{eq:u_indep}:
\begin{align*}
0\le u\le 1 \quad \mbox{a.e.}
\end{align*}
\end{lemma}

\proof
We choose $\varphi=-u_- \in L^2(0,T; H^1(\Omega))$ as a test function in the weak form of the equation, with $T$ replaced by an arbitrary $t \in (0,T]$. 
Here, 
we use the notation $u_{\pm}:= \max \{ \pm u, 0\}$. Then, we have 
\[
\int_0^t \la \frac{\partial u}{\partial t}, -u_- \ra 
+\int_0^t \la   \beta_\eta'(u) \nabla u ,\nabla (-u_-)\ra
 -  \int_0^t \la g(u) \nabla \sum_{j=1}^M a_j v_j,\nabla (-u_-)\ra
 =0.  
\]
{
Since $\beta_\eta'(u) \geq \eta$ and $\nabla(-u_-)=\boldsymbol{1}_{\{u<0\}}\nabla u$, the second term on the left-hand side satisfies
\[
\int_0^t \la \beta_\eta'(u)\nabla u,\nabla(-u_-)\ra
=\int_{\{u<0\}\cap Q_t}\beta_\eta'(u)|\nabla u|^2\geq0.
\]
Moreover, $g(u)=0$ on $\{u<0\}$ by \eqref{assum_g}, whereas $\nabla(-u_-)=0$ on $\{u\geq0\}$; hence, the third term vanishes. By the chain rule and $u_{0-}=0$, which follows from \eqref{assum_init_u}, we obtain
}
\begin{align*}
\int_0^t \la \frac{\partial u}{\partial t}, -u_- \ra &=
\int_0^t \la \frac{\partial (-u_-)}{\partial t}, -u_- \ra 
=
\frac{1}{2} \ld u_-(t) \rd^2 - \frac{1}{2} \ld u_{0-} \rd^2 
=
\frac{1}{2} \ld u_-(t) \rd^2 \le 0. 
\end{align*}
Thus, we get $u_-=0$ a.e. and then $u\ge 0$ a.e. 

The weak form of \eqref{eq:u_indep} can be rewritten as follows: 
\[
\int_0^t \la \frac{\partial }{\partial t}(u-1), \varphi \ra 
+\int_0^t \la \beta_\eta' (u) \nabla (u-1),\nabla \varphi\ra
 -  \int_0^t \la g(u) \nabla \sum_{j=1}^M a_j v_j,\nabla \varphi\ra
 =0 
\]
for arbitrary $t \in (0,T]$ and $\varphi\in L^2(0,T;H^1(\Omega))$. 
We set $\varphi = (u-1)_+$. Because $g(u)=0$ where $u>1$ and $\nabla (u-1)_+ =0$ where $u\le 1$, we similarly obtain the following:
\[
\frac{1}{2} \ld (u-1)_+(t) \rd^2 \le \frac{1}{2} \ld (u_{0}-1)_+ \rd^2=0. 
\]
{ Here, $(u_0-1)_+=0$ follows from \eqref{assum_init_u}.}
Thereby, we obtain $u\le 1$. {Since $\widetilde\beta=\beta$ on $[0,1]$, the auxiliary modification does not affect the equation once the invariant-region estimate $0\le u\le1$ has been established. Accordingly, all subsequent estimates and the limit $\eta\to0$ may be written in terms of the original function $\beta$, and no additional assumption beyond Assumption~\eqref{assum_gamma} is required.} This completes the proof. 
\endproof

\begin{lemma}
\label{Lem_estim_vj}
{If $v_j$ is the solution of \eqref{eq:vj_indep}, then} there exists a positive constant $C$ independent of $\xi$, $\eta$ and $\bar{u}$ but dependent on $|\Omega|, T, \{ d_j\}$ and { $\ld v_{0j}\rd_{H^1(\Omega)}$} such that 
\begin{equation*}
\begin{aligned}
	\sqrt{\xi} \ld v_j \rd_{L^\infty(0,T; H^1(\Omega))} +
    \ld v_{j}\rd_{L^2(0,T; H^1(\Omega))}
    +
    \ld \Delta v_{j}\rd_{L^2(Q_T)}
     \le    C. 
\end{aligned}
\end{equation*}
\end{lemma}

Throughout this section, $C$ denotes a generic positive constant independent of $\eta$, $\xi$ and $\bar{u}$. { It may depend on $|\Omega|$, $T$, $M$, $\{a_j\}_{j=1}^M$, $\{d_j\}_{j=1}^M$, $L_g$, $\beta$, and the initial data $\{v_{0j}\}_{j=1}^M$.}

\proof
Multiplying the equation for $v_j$ by $v_j$, integrating both sides in space, and using integration by parts and the Cauchy--Schwarz inequality, we have 
\begin{equation*}
\begin{aligned}
	\frac{\xi}{2} \frac{ d }{ dt }\ld v_j\rd^2
    +
    d_j\ld \nabla v_{j}\rd^2
    + 
    \ld  v_{j}\rd^2
     =
    \la \bar{u}, v_j \ra \\ 
     \le    
    \frac{1}{2}\ld  \bar{u}\rd^2 
    +
    \frac{1}{2}\ld  v_{j}\rd^2
     \le    
    \frac{1}{2}|\Omega| 
    +
    \frac{1}{2}\ld  v_{j}\rd^2.    
\end{aligned}
\end{equation*}
Absorbing the last term on the right-hand side into the left-hand side and integrating both sides with respect to time over $[0,t]$ for $t \in (0,T]$, we obtain the following. 
\begin{equation*}
\begin{aligned}
	\frac{\xi}{2} \ld v_j(t) \rd^2
    +
    d_j\ld \nabla v_{j}\rd_{L^2(Q_t)}^2
    + 
    \frac{1}{2}\ld  v_{j}\rd_{L^2(Q_t)}^2
     \le    
    \frac{1}{2}|\Omega|T 
    +
    \frac{\xi}{2} \ld v_{0j} \rd^2. 
\end{aligned}
\end{equation*}
Therefore, we get 
\begin{equation*}
\begin{aligned}
	\sqrt{\xi} \ld v_j \rd_{L^\infty(0,T; L^2(\Omega))} 
    +
    \ld v_{j}\rd_{L^2(0,T; H^1(\Omega))}
     \le    C. 
\end{aligned}
\end{equation*}

Multiply the equation for $v_j$ by ${-}\Delta 
v_j$ and apply a strategy similar to the above to get 
\begin{equation*}
\begin{aligned}
	\frac{\xi}{2} \ld \nabla v_j(t) \rd^2
    +
    \frac{d_j}{2} \ld \Delta v_{j}\rd_{L^2(Q_t)}^2
    + 
    \ld  \nabla v_{j}\rd_{L^2(Q_t)}^2
     \le    
    \frac{1}{2d_j}|\Omega|T 
    +
    \frac{\xi}{2} \ld \nabla v_{0j} \rd^2. 
\end{aligned}
\end{equation*}
Hence, we have 
\begin{equation*}
\begin{aligned}
	\sqrt{\xi} \ld v_j \rd_{L^\infty(0,T; H^1(\Omega))} 
    +
    \ld \Delta v_{j}\rd_{L^2(Q_T)}
     \le    C. 
\end{aligned}
\end{equation*}
Thus, we complete the proof. 
\endproof

\begin{lemma}
\label{Lem_estim_vj_H1L2}
There exists a positive constant $C_\xi$ independent of $\eta$ and $\bar{u}$ but dependent on $\xi$, $|\Omega|, T, \{ d_j\}$ and { $\ld v_{0j}\rd_{H^1(\Omega)}$} such that 
\begin{equation*}
\begin{aligned}
	 \ld v_j \rd_{H^1(0,T; L^2(\Omega))}
     \le    C_\xi. 
\end{aligned}
\end{equation*}
\end{lemma}

\proof
Multiply the equation for $v_j$ by $\frac{\partial v_j}{\partial t}$, integrate both sides in $Q_t$ for $t\in (0,T]$, and use integration by parts to obtain  
\begin{align}
   \xi \ld \frac{\partial v_j}{\partial t}\rd_{L^2(Q_t)}^2
    +
    \frac{d_j}{2} \int_0^t \frac{d}{dt}\ld \nabla v_{j}\rd^2
    + 
    \frac{1}{2}\int_0^t \frac{d}{dt}\ld v_{j}\rd^2
     =
    \int_0^t \la \bar{u}, \frac{\partial v_j}{\partial t} \ra \label{estim:vj:H1L2:xi}\\
    \le 
    \frac{1}{2\xi} \ld \bar{u} \rd_{L^2(Q_t)}^2
    +
    \frac{\xi}{2} \ld \frac{\partial v_j}{\partial t}\rd_{L^2(Q_t)}^2. \notag
\end{align}
Thus, we have 
\begin{equation*}
\begin{aligned}
   \frac{\xi}{2} \ld \frac{\partial v_j}{\partial t}\rd_{L^2(Q_t)}^2
    +
    \frac{d_j}{2} \ld \nabla v_{j}(t) \rd^2
    + 
    \frac{1}{2}\ld v_{j}(t)\rd^2
    \le 
    \frac{d_j}{2} \ld \nabla v_{0j} \rd^2
    + 
    \frac{1}{2}\ld v_{0j}\rd^2
    +
    \frac{1}{2\xi} |\Omega|T, 
\end{aligned}
\end{equation*}
that completes the proof.
\endproof

\begin{lemma}
\label{Lem_estim_u_L2H1}
{If $u$ is the solution of \eqref{eq:u_indep}, then} the following estimate holds: 
\begin{equation*}
\begin{aligned}
\ld  \beta(u)\rd_{L^2(0,T;H^1(\Omega))}
+
\ld  \beta_\eta(u)\rd_{L^2(0,T;H^1(\Omega))}
+
\ld  u\rd_{H^1(0,T;H^1(\Omega)^*)}
\le
C.
\end{aligned}
\end{equation*}
\end{lemma}

\proof
{Choosing $\varphi=\beta(u) \in L^2(0,T; H^1(\Omega))$ as a test function in the weak form of the equation, we obtain the following relation.}
\[
\int_0^T \la \frac{\partial u}{\partial t}, \beta(u) \ra 
+\int_0^T \la   \nabla \lr \eta u+\beta(u) \rr ,\nabla \beta(u) \ra
 -  \int_0^T \la g(u) \nabla \sum_{j=1}^M a_j v_j,\nabla \beta(u)\ra
 =0.  
\]
Define 
\[
\Phi (s):= \int_0^s \beta(r)\, dr \quad s\in \R. 
\]
It is easy to see that $\Phi$ is convex and satisfies 
\begin{equation}
\label{prop:Phi}
0\le \frac{1}{2L_\beta}\beta(s)^2 \le \Phi (s) \le \frac{L_\beta}{2}s^2 \qquad \mbox{for } s\in [0,1], 
\end{equation}
where $L_\beta:= \max_{s\in [0,1]} \beta'(s)$. 
{
Since
$
\la \frac{\partial u}{\partial t}, \beta(u) \ra
=
\frac{d}{dt}\int_\Omega \Phi(u),
$
and, by the monotonicity of $\beta$,
\begin{align*}
\int_0^T \la \nabla\lr \eta u+\beta(u)\rr,\nabla\beta(u)\ra\,dt
&=
\ld\nabla\beta(u)\rd_{L^2(Q_T)}^2
+
\eta\int_0^T\!\!\int_\Omega \beta'(u)|\nabla u|^2\,dx\,dt\\
&\ge
\ld\nabla\beta(u)\rd_{L^2(Q_T)}^2,
\end{align*}
holds, the Cauchy--Schwarz inequality yields
}
\[
\ld  \nabla \beta(u)\rd_{L^2(Q_T)}^2
\le
\int_\Omega \Phi(u_0)
+
\frac{1}{2} \ld g(u) \nabla \sum_{j=1}^M a_j v_j \rd_{L^2(Q_T)}^2 
+
\frac{1}{2} \ld \nabla \beta(u) \rd_{L^2(Q_T)}^2. 
\]
Therefore, it follows from \eqref{prop:Phi}, the continuity of $g$ and Lemma~\ref{Lem_estim_vj} that 
\[
\ld  \nabla \beta(u)\rd_{L^2(Q_T)}
\le
C.
\]
Similarly, by choosing $\beta_\eta(u)$ as a test function, we obtain the following.
\begin{equation}
\label{estim:beta_eta}
\ld  \nabla \beta_\eta(u)\rd_{L^2(Q_T)}
\le
C.
\end{equation}

From the weak form of the equation, we have the following estimate for arbitrary $\psi \in H^1(\Omega)$ and for a.e. $t\in (0,T]$:
\begin{align*}
\left| \la \frac{\partial u}{\partial t}(t), \psi \ra \right| 
&\le 
\left| \la \nabla \beta_\eta(u(t)),\nabla \psi\ra \right|
 +
   \left| \la g(u(t)) \nabla \sum_{j=1}^M a_j v_j(t),\nabla \psi\ra \right| \\
   &\le 
 \ld \nabla \beta_\eta(u(t)) \rd \ld \nabla \psi\rd 
 +
   L_g \ld \nabla \sum_{j=1}^M a_j v_j(t) \rd \ld \nabla \psi\rd. 
\end{align*}
Thus, we get
\begin{align*}
\ld \frac{\partial u}{\partial t}(t) \rd_{H^1(\Omega)^*}^2  &=\lr \sup_{\ld \psi\rd_{H^1(\Omega)}\neq 0} 
\frac{\left| \la \frac{\partial u}{\partial t}(t), \psi \ra \right|}{\ld \psi\rd_{H^1(\Omega)}}\rr^2\\
&
\le 
 2\ld \nabla \beta_\eta(u(t)) \rd^2  +
   2 L_g^2 \ld \nabla \sum_{j=1}^M a_j v_j(t) \rd^2. 
\end{align*}
Integrate both sides in $t$ from $0$ to $T$ and use \eqref{estim:beta_eta} and Lemma~\ref{Lem_estim_vj} to obtain the desired $H^1(0,T;H^1(\Omega)^*)$-estimate. 
\endproof

When $\eta = 0$, the diffusion is degenerate in general, so obtaining a uniform estimate for $u$ with respect to $\eta$ in $L^2(0,T; H^1(\Omega))$ cannot be expected.
Therefore, we have the following result.

\begin{lemma}
\label{Lem_estim_u_L2H1_dep_eta}
{If $u$ is the solution of \eqref{eq:u_indep}, then} there exists a positive constant $C_\eta$ that depends on $\eta$ such that
\begin{equation*}
\begin{aligned}
\ld  u\rd_{L^2(0,T;H^1(\Omega))}
\le
C_\eta.
\end{aligned}
\end{equation*}
\end{lemma}

{
\proof
For each fixed $\eta\in(0,1]$, the desired estimate follows immediately from Lemma~\ref{Lem_estim_u_L2H1} and the identity $\eta u=\beta_\eta(u)-\beta(u)$.
\endproof
}

We are now ready to prove the existence of a weak solution to the non-degenerate problem~\eqref{eq:KSrelax}.

\begin{theorem}
\label{th_exist_KSrelax}
Under Assumptions \eqref{assum_gamma}--\eqref{assum_init_v}, there exists a weak solution of \eqref{eq:KSrelax}. 
\end{theorem}
\proof
Consider a mapping $\Theta(\bar{u}) = u$, where $u$ is the weak solution to \eqref{eq:u_indep} obtained via \eqref{eq:vj_indep} for a given $\bar{u} \in \mathcal{K}$. Lemma~\ref{Lem_linf_u} implies that $\Theta$ maps $\mathcal{K}$ into itself. 
{In view of Lemmas~\ref{Lem_estim_u_L2H1} and~\ref{Lem_estim_u_L2H1_dep_eta}, $\Theta(\mathcal{K})$ is uniformly bounded in $L^2(0,T;H^1(\Omega))\cap H^1(0,T;H^1(\Omega)^*)$. Hence, the Aubin--Lions lemma implies that $\Theta(\mathcal{K})$ is relatively compact in $L^2(Q_T)$.
We next prove that $\Theta$ is continuous.} To this end, let $\{ \bar{u}^n \}_{n=1}^\infty$ be a sequence of functions in $\mathcal{K}$ and $\bar{u}\in \mathcal{K}$ be such that $\bar{u}^n \to \bar{u}$ in $L^2(Q_T)$ as $n\to \infty$. 
For each $n$, let $v_j^n$ ($j=1,2,\dots,M$) be the solution of \eqref{eq:vj_indep} with $\bar{u}$ replaced by $\bar{u}^n$, and let $u^n$ denote the weak solution of \eqref{eq:u_indep} using these $v_j^n$. 
In view of Lemmas~\ref{Lem_estim_vj}-\ref{Lem_estim_u_L2H1_dep_eta}, $u^n$ and { all $v_j^n$, $j=1,\ldots,M$,} are uniformly bounded in $L^2(0,T;H^1(\Omega))$ and $H^1(0,T;H^1(\Omega)^*)$ with respect to $n$. Moreover, $\beta_\eta(u^n)$ is uniformly bounded in $L^2(0,T;H^1(\Omega))$ with respect to $n$. 
From the Aubin--Lions Lemma, the injection of $L^2(0,T;H^1(\Omega))\cap H^1(0,T;H^1(\Omega)^*)$ into $L^2(Q_T)$ is compact. 
{ Since $M$ is fixed and finite, a finite diagonal extraction yields a common subsequence for $u^n$ and all $v_j^n$, $j=1,\ldots,M$.} Hence, there exist { a subsequence, which is denoted by $\{u^{n}\}$ and $\{v_j^{n}\}$ again,} and functions $u^*\in \mathcal{K}\cap L^2(0,T;H^1(\Omega))\cap H^1(0,T;H^1(\Omega)^*)$ and $v_j^*\in L^\infty(0,T;H^1(\Omega))\cap H^1(0,T;L^2(\Omega))$ such that 
\begin{align*}
&u^{n} \to u^*, \quad v_j^{n} \to v_j^* 
& & \mbox{strongly in } L^2(Q_T), \mbox{ a.e. in } Q_T,   \notag \\
& & &  \mbox{weakly in } L^2(0,T;H^1(\Omega)) \mbox{ and } H^1(0,T;H^1(\Omega)^*), \\
&\beta_\eta (u^{n}) \to \beta_\eta(u^*) 
& & \mbox{strongly in } L^2(Q_T), \mbox{ a.e. in } Q_T,   \notag \\
& & &\mbox{weakly in } L^2(0,T;H^1(\Omega)) \notag 
\end{align*}
as $n$ tends to infinity. 
Since $g$ is continuous, $g(u^{n})$ converges to $g(u^*)$ in $L^2(Q_T)$ and a.e. as $n$ tends to infinity. 
In the weak form of the equation for $v_j^n$, i.e., 
\begin{align}
\xi \int_0^T \la \frac{\partial v_j^n}{\partial t}, \varphi \ra 
+d_j\int_0^T \la   \nabla v_j^n  ,\nabla \varphi \ra
 +  \int_0^T \la v_j^n, \varphi\ra
 = \int_0^T \la \bar{u}^n, \varphi\ra \label{weak_vjn}  
\end{align}
for $\varphi \in L^2(0,T;H^1(\Omega))$, {passing to the limit as $n\to\infty$, we see that} $v_j^*$ is the weak solution of \eqref{eq:vj_indep} satisfying 
\begin{align}
\xi \int_0^T \la \frac{\partial v_j^*}{\partial t}, \varphi \ra 
+d_j\int_0^T \la   \nabla v_j^*  ,\nabla \varphi \ra
 +  \int_0^T \la v_j^*, \varphi\ra
 = \int_0^T \la \bar{u}, \varphi\ra.  \label{weak_vj*}
\end{align}
From the standard regularity theory~\citep{brezis}, $v_j^*\in L^2(0,T; H^2(\Omega))$, and thus, $v_j^*$ is the solution of \eqref{eq:vj_indep}. 
Subtracting \eqref{weak_vj*} from \eqref{weak_vjn} and choosing $\varphi = v_j^n-v_j^*$, { since $v_j^n(0)=v_j^*(0)=v_{0j}$, the initial-time term vanishes,} we obtain 
\begin{equation*}
\begin{aligned}
	\frac{\xi}{2} \ld v_j^n(T)-v_j^*(T) \rd^2
    +
    d_j\ld \nabla v_j^n- \nabla v_j^*\rd_{L^2(Q_T)}^2
    + 
    \ld  v_j^n-v_j^* \rd_{L^2(Q_T)}^2\\
     =
    \int_0^T \la \bar{u}^n - \bar{u}, v_j^n-v_j^* \ra 
    \to 0 \mbox{ as } n\to \infty,     
\end{aligned}
\end{equation*}
which implies the strong convergence of $v_j^n$ in $L^2(0,T;H^1(\Omega))$, and this further leads to a.e. convergence. 
{
Moreover,
\begin{align*}
&\left\|
g(u^n)\nabla\sum_{j=1}^M a_jv_j^n
-g(u^*)\nabla\sum_{j=1}^M a_jv_j^*
\right\|_{L^2(Q_T)}\\
&\quad\le
\|g(u^n)\|_{L^\infty(Q_T)}
\left\|\nabla\sum_{j=1}^M a_j(v_j^n-v_j^*)\right\|_{L^2(Q_T)}
+
\left\|
\bigl(g(u^n)-g(u^*)\bigr)
\nabla\sum_{j=1}^M a_jv_j^*
\right\|_{L^2(Q_T)}\\
&\quad \to 0.
\end{align*}
Indeed, the first term tends to zero by the strong convergence of $v_j^n$ in $L^2(0,T;H^1(\Omega))$ and the boundedness of $g$ on $[0,1]$. The second term tends to zero by the pointwise convergence of $g(u^n)$, the boundedness of $g$, and the dominated convergence theorem.
}
Therefore, taking the limit as $n \to \infty$ in the weak form of the equation for $u^n$, we conclude that $u^*$ is a weak solution to \eqref{eq:u_indep}.
{By the uniqueness proved above, $u^*=\Theta(\bar u)$; since every subsequence admits a further subsequence converging to this same limit, the whole sequence converges to $\Theta(\bar u)$.
Consequently, $\Theta$ is a continuous map on $\mathcal{K}$.} 
Now, we can apply Schauder's fixed point theorem and find that $\Theta$ has a fixed point, that is, there exists $u$ such that $\Theta(u)=u$, which indicates the existence of a weak solution to \eqref{eq:KSrelax}. 
\endproof

\subsection{Non-degenerate to degenerate limit}

In the previous subsection, the existence of a weak solution $\lr u^{\eta,\xi,M}, v_1^{\eta,\xi,M}, \dots, v_M^{\eta,\xi,M} \rr$ of the non-degenerate problem~\eqref{eq:KSrelax} was established. In this subsection, we demonstrate the existence of a weak solution of \eqref{eq:KSM} by considering the limit as the regularization parameter $\eta$ tends to zero. 
To this end, $\eta$-independent a priori estimates are required, but thanks to Lemmas~\ref{Lem_linf_u}, \ref{Lem_estim_vj}, \ref{Lem_estim_vj_H1L2}, and \ref{Lem_estim_u_L2H1}, we have already obtained the following estimates.

\begin{lemma}
\label{lem:eta_uniform_estimates}
The following estimates hold for the weak solution $\lr u^{\eta,\xi,M}, v_1^{\eta,\xi,M}, \dots, v_M^{\eta,\xi,M} \rr$ of \eqref{eq:KSrelax}:
\begin{gather*}
0\le u^{\eta,\xi,M}\le 1 \quad \mbox{a.e.}, \\
\ld  \beta \lr u^{\eta,\xi,M} \rr \rd_{L^2(0,T;H^1(\Omega))}
+
\ld  u^{\eta,\xi,M}\rd_{H^1(0,T;H^1(\Omega)^*)}
\le
C, \\
	\sqrt{\xi} \ld v_j^{\eta,\xi,M} \rd_{L^\infty(0,T; H^1(\Omega))} +
    \ld v_{j}^{\eta,\xi,M}\rd_{L^2(0,T; H^1(\Omega))}
    +
    \ld \Delta v_{j}^{\eta,\xi,M}\rd_{L^2(Q_T)}
     \le    C, \\
	 \ld v_j^{\eta,\xi,M} \rd_{H^1(0,T; L^2(\Omega))}
     \le    C_\xi. 
\end{gather*}
Here, $C$ is a positive constant independent of $\xi$ and $\eta$ but dependent on $|\Omega|$, $T$, {$M$, $\{a_j\}_{j=1}^M$}, $\{d_j\}_{j=1}^M$, {$L_g$, $\beta$}, and { $\ld v_{0j}\rd_{H^1(\Omega)}$}. 
The positive constant $C_\xi$ depends on $\xi$ but is independent of $\eta$. 
\end{lemma}

{
The following lemma is a special case of the compactness argument of
\citet[Lemma~1.9]{altLuckhaus1983}. We include a short proof adapted
to the present periodic setting.

\begin{lemma}[Alt--Luckhaus-type compactness]
{ Assume that $\beta\in C^1([0,1])$ is  a strictly increasing function.}
\label{lem:nonlinear_compactness}
Let $\{z_n\}$ be a sequence satisfying
\[
0\le z_n\le 1 \quad\mbox{a.e. in }Q_T
\]
and
\[
\sup_n\left(
 \left\|z_n\right\|_{H^1(0,T;H^1(\Omega)^*)}
 +
 \left\|\beta(z_n)\right\|_{L^2(0,T;H^1(\Omega))}
\right)<+\infty.
\]
Then $\{\beta(z_n)\}$ is relatively compact in $L^p(Q_T)$
for every $1\le p<+\infty$.
\end{lemma}

\begin{proof}
For $0<h<T$, set
\[
\delta_h z_n(t):=z_n(t+h)-z_n(t).
\]
Since
\[
\delta_h z_n(t)
=
\int_t^{t+h}\frac{\partial z_n}{\partial s}(s)\,ds,
\]
we have
\[
\left\|\delta_h z_n\right\|_
 {L^2(0,T-h;H^1(\Omega)^*)}
\le
h\left\|\frac{\partial z_n}{\partial t}\right\|_
 {L^2(0,T;H^1(\Omega)^*)}.
\]
Moreover,
\[
\left\|\delta_h\beta(z_n)\right\|_
 {L^2(0,T-h;H^1(\Omega))}
\le
2\left\|\beta(z_n)\right\|_
 {L^2(0,T;H^1(\Omega))}.
\]
Therefore,
\begin{align*}
&\int_0^{T-h}\!\!\int_\Omega
 \delta_h z_n(t)\,
 \delta_h\beta(z_n)(t)\,dx\,dt
\\
&\qquad\le
\left\|\delta_h z_n\right\|_
 {L^2(0,T-h;H^1(\Omega)^*)}
\left\|\delta_h\beta(z_n)\right\|_
 {L^2(0,T-h;H^1(\Omega))}
\le Ch.
\end{align*}
Let
\[
L_\beta:=\max_{s\in[0,1]}\beta'(s).
\]
Since $\beta$ is increasing and Lipschitz continuous on $[0,1]$,
\[
|\beta(r)-\beta(s)|^2
\le
L_\beta(r-s)\bigl(\beta(r)-\beta(s)\bigr)
\qquad (r,s\in[0,1]).
\]
It follows that
\[
\left\|\delta_h\beta(z_n)\right\|_{L^2(Q_{T-h})}^2
\le Ch.
\]
On the other hand, the uniform boundedness of $\beta(z_n)$ in
$L^2(0,T;H^1(\Omega))$ gives, using the periodicity of $\Omega$,
\[
\left\|
 \beta(z_n)(\,\cdot+y,\cdot\,)-\beta(z_n)
\right\|_{L^2(Q_T)}
\le
|y|\left\|\nabla\beta(z_n)\right\|_{L^2(Q_T)}.
\]
Thus, both the space and time translates of $\beta(z_n)$ tend to zero
uniformly with respect to $n$. In addition, $0\le z_n\le1$ implies a
uniform $L^\infty(Q_T)$ bound for $\beta(z_n)$, which controls the time
boundary strips. The Riesz--Fr\'echet--Kolmogorov compactness criterion
(see, e.g., \citet{Simon1987})
yields the relative compactness of $\{\beta(z_n)\}$ in $L^2(Q_T)$.
Since $\{\beta(z_n)\}$ is uniformly bounded in $L^\infty(Q_T)$,
an interpolation argument yields its relative compactness in $L^p(Q_T)$
for every $1\le p<+\infty$.
\end{proof}
}

\begin{proof}[Proof of Theorem~\ref{thm:exist:KSM}]
{
Set $\widehat\eta_n=1/n$. Applying
Lemma~\ref{lem:nonlinear_compactness} to
$\{u^{\widehat\eta_n,\xi,M}\}_{n\in\N}$, we conclude that
$\{\beta(u^{\widehat\eta_n,\xi,M})\}_{n\in\N}$ is relatively compact in
$L^2(Q_T)$. Hence, using the estimates above and taking a common
subsequence for the finite family of chemical fields, there exist a
subsequence $\{\eta_k\}$ of $\{\widehat\eta_n\}$,
$u^*\in\mathcal{K}\cap H^1(0,T;H^1(\Omega)^*)$,
$\beta^*\in L^2(0,T;H^1(\Omega))$, and
$v_j^*\in L^\infty(0,T;H^1(\Omega))
\cap H^1(0,T;L^2(\Omega))
\cap L^2(0,T;H^2(\Omega))$ such that
\begin{align}
&u^{\eta_k,\xi,M} \rightharpoonup u^* 
& &  \mbox{weakly in } L^2(Q_T) \mbox{ and } H^1(0,T;H^1(\Omega)^*), \label{conv:u_eta} \\
&\beta \lr u^{\eta_k,\xi,M}\rr \to \beta^* 
& & \mbox{strongly in } L^2(Q_T), \mbox{ a.e. in } Q_T,   \label{conv:beta_eta} \\
& & & \mbox{weakly in } L^2(0,T;H^1(\Omega)), \notag \\
&v_j^{\eta_k,\xi,M} \to v_j^* 
& & \mbox{strongly in } L^2(0,T;H^1(\Omega)), \mbox{ a.e. in } Q_T,   \notag \\
&v_j^{\eta_k,\xi,M} \stackrel{*}{\rightharpoonup} v_j^*
& & \mbox{weakly-* in }L^\infty(0,T;H^1(\Omega)), \notag\\
&v_j^{\eta_k,\xi,M} \rightharpoonup v_j^*
& & \mbox{weakly in }H^1(0,T;L^2(\Omega))
\mbox{ and }L^2(0,T;H^2(\Omega)) \notag
\end{align}
}
as $\eta_k$ tends to zero. 
{
Together with the weak convergence of $u^{\eta_k,\xi,M}$,
Lemma~6.1 in \citet{eymard1999fvm} gives
$\beta^*=\beta(u^*)$.
Since $\beta$ is strictly increasing, its inverse is continuous on
$\beta([0,1])$. Consequently,
\[
u^{\eta_k,\xi,M}\to u^*
\quad\mbox{a.e. in }Q_T.
\]
The uniform bound $0\le u^{\eta_k,\xi,M}\le1$ then implies
\[
u^{\eta_k,\xi,M}\to u^*
\quad\mbox{strongly in }L^2(Q_T).
\]
}
{We also obtain the convergence} of $g\lr u^{\eta_k,\xi,M}\rr \nabla \sum_{j=1}^M a_j v_j^{\eta_k,\xi,M} 
\to 
g(u^*) \nabla \sum_{j=1}^M a_j v_j^*$ in $L^2(Q_T)$ { as $\eta_k$ tends to zero}. 
Now, we can pass to the limit in $\eta_k$ in the following weak formulation:
\begin{align*}
\int_0^T \la \frac{\partial u^{\eta_k,\xi,M}}{\partial t}, \varphi \ra 
&
-\int_0^T \la \eta_k u^{\eta_k,\xi,M}, \Delta \varphi\ra
+\int_0^T \la \nabla \beta \lr u^{\eta_k,\xi,M} \rr,\nabla \varphi\ra \\
 & \hspace*{2cm}
 -  \int_0^T \la g\lr u^{\eta_k,\xi,M} \rr \nabla \sum_{j=1}^M a_j v_j^{\eta_k,\xi,M},\nabla \varphi\ra
 =0, \\ 
\xi \int_0^T \la \frac{\partial v_j^{\eta_k,\xi,M}}{\partial t}, \psi \ra 
&+d_j\int_0^T \la   \nabla v_j^{\eta_k,\xi,M}  ,\nabla \psi \ra
 +  \int_0^T \la v_j^{\eta_k,\xi,M}, \psi \ra
 = \int_0^T \la u^{\eta_k,\xi,M}, \psi \ra 
\end{align*}
for all functions $\varphi\in L^2(0,T;H^2(\Omega))$ and $\psi\in L^2(0,T;H^1(\Omega))$. 
{ Each term in the limiting weak formulation for $u^*$ is continuous with respect to the $L^2(0,T;H^1(\Omega))$ norm of the test function. Since $H^2(\Omega)$ is dense in $H^1(\Omega)$ under the periodic boundary condition, the identity extends to every $\varphi\in L^2(0,T;H^1(\Omega))$.}
{
The weak convergences in $H^1(0,T;H^1(\Omega)^*)$ and
$H^1(0,T;L^2(\Omega))$, together with the continuity of the trace
operators, yield
\[
u^*(0)=u_0,
\qquad
v_j^*(0)=v_{0j}.
\]
}
From the regularities of $u^*$ and $v_j^*$, we realize that $\lr u^*, v_1^*, \dots, v_M^* \rr$ is the weak solution of \eqref{eq:KSM}. 
\end{proof}


\section{Parabolic--parabolic to parabolic--elliptic limit of the Keller--Segel type system}
\label{sec:pptope}
In this section, we consider the relationship between \eqref{eq:KSM} and \eqref{eq:KSPE}, particularly focusing on the limit as the relaxation parameter $\xi$ tends to zero.

{
The fast signal diffusion limit from parabolic--parabolic to parabolic--elliptic Keller--Segel systems has been studied for systems with linear or non-degenerate nonlinear diffusion; see, e.g., \citep{Mizukami2019,Freitag2020,OgawaSuguro2023}. \citet{Mizukami2019} considers the standard Keller--Segel system with linear diffusion on a bounded domain, \citet{Freitag2020} treats a system with non-degenerate nonlinear diffusion, and \citet{OgawaSuguro2023} establish strong convergence in uniformly local spaces on the whole space. Here, we establish this limit in a weak-solution framework for a finite family of chemical fields under assumptions allowing degenerate nonlinear diffusion and density saturation.
}

We consider the weak solution $\lr u^{\xi,M}, v_1^{\xi,M}, \dots, v_M^{\xi,M} \rr$ of the {parabolic--parabolic} problem~\eqref{eq:KSM}. 
{
For any such weak solution, the estimates for $v_j^{\xi,M}$ follow
directly from its equation, exactly as in Lemma~3.2. Moreover,
Definition~2.2 ensures that $\beta(u^{\xi,M})$ is an admissible test
function, and the nonlinear chain rule of
\citet[Lemma~1.5]{altLuckhaus1983} justifies
the corresponding time term. Hence, the estimates for $u^{\xi,M}$
follow as in Lemma~3.4.
}
\begin{lemma}
\label{lem:KSM:arpriori}
The following estimates hold for {every} weak solution $\lr u^{\xi,M}, v_1^{\xi,M}, \dots, v_M^{\xi,M} \rr$ of \eqref{eq:KSM}:
\begin{gather*}
0\le u^{\xi,M}\le 1 \quad \mbox{a.e.}, \\
{\ld  \beta \lr u^{\xi,M} \rr \rd_{L^2(0,T;H^1(\Omega))}}
+
\ld  u^{\xi,M}\rd_{H^1(0,T;H^1(\Omega)^*)}
\le
C, \\
\sqrt{\xi} \ld v_j^{\xi,M} \rd_{L^\infty(0,T; H^1(\Omega))} +
    \ld v_{j}^{\xi,M}\rd_{L^2(0,T; H^1(\Omega))}
    +
    \ld \Delta v_{j}^{\xi,M}\rd_{L^2(Q_T)}
     \le    C. 
\end{gather*}
Here, $C$ is a positive constant independent of $\xi$ but dependent on $|\Omega|$, $T$, {$M$, $\{a_j\}_{j=1}^M$}, $\{d_j\}_{j=1}^M$, {$L_g$, $\beta$}, and { $\ld v_{0j}\rd_{H^1(\Omega)}$}. 
\end{lemma}

In addition, we establish the following uniform boundedness.
\begin{lemma}
\label{lem:KSM:arpriori_v}
There exists a positive constant $C$ independent of $\xi$ such that the following estimate holds. 
\begin{gather*}
\xi^{3/4} \ld \frac{\partial v_j^{\xi,M}}{\partial t}\rd_{L^2(Q_T)}
+
\xi^{1/4} \ld v_j^{\xi,M} \rd_{L^\infty(0,T; H^1(\Omega))}
     \le    C. 
\end{gather*}
\end{lemma}

\begin{proof}
Apply the same strategy as in \eqref{estim:vj:H1L2:xi} and use integration by parts and Lemma~\ref{lem:KSM:arpriori}.
{The integration by parts in time is justified by $u^{\xi,M}\in H^1(0,T;H^1(\Omega)^*)$ and $v_j^{\xi,M}\in L^2(0,T;H^1(\Omega))\cap H^1(0,T;L^2(\Omega))$.}
We obtain
\begin{align*}
   \xi& \ld \frac{\partial v_j^{\xi,M}}{\partial t}\rd_{L^2(Q_t)}^2
    +
    \frac{d_j}{2} \int_0^t \frac{d}{dt}\ld \nabla v_{j}^{\xi,M}\rd^2
    + 
    \frac{1}{2}\int_0^t \frac{d}{dt}\ld v_{j}^{\xi,M}\rd^2
     =
    \int_0^t \la u^{\xi,M}, \frac{\partial v_j^{\xi,M}}{\partial t} \ra \\
&= 
    \la u^{\xi,M}(t),  v_j^{\xi,M}(t) \ra
    -
    \la u_0,  v_{0j} \ra
{ -}
\int_0^t \la  \frac{\partial u^{\xi,M}}{\partial t}, v_j^{\xi,M} \ra \\
&\le 
    \ld u^{\xi,M}(t)\rd  \ld v_j^{\xi,M}(t) \rd
    +
    \ld u_0 \rd \ld v_{0j} \rd
+
{\ld  \frac{\partial u^{\xi,M}}{\partial t}\rd_{L^2(0,T; H^1(\Omega)^*)}}  \ld v_j^{\xi,M} \rd_{L^2(0,T; H^1(\Omega))}\\
&\le 
    { |\Omega|^{1/2}}  \ld v_j^{\xi,M}(t) \rd
    +
C'.
\end{align*}
Here, $C'$ is a positive constant independent of $\xi$ that bounds the last two terms.
{ Multiplying both sides by $\sqrt{\xi}$}, we have 
\begin{align*}
   \xi^{3/2}& \ld \frac{\partial v_j^{\xi,M}}{\partial t}\rd_{L^2(Q_t)}^2
    +
    \frac{d_j}{2} \sqrt{\xi} \ld \nabla v_{j}^{\xi,M}(t)\rd^2
    + 
    \frac{1}{2}\sqrt{\xi}\ld v_{j}^{\xi,M}(t)\rd^2\\
&\le 
    \frac{d_j}{2} \sqrt{\xi} \ld \nabla v_{0j}\rd^2
+
    \frac{1}{2}\sqrt{\xi} \ld v_{0j}\rd^2
    +
    { |\Omega|^{1/2}} \sqrt{\xi}  \ld v_j^{\xi,M}(t) \rd
    +
C'\sqrt{\xi}.
\end{align*}
{ Since $
\sqrt{\xi}\ld v_j^{\xi,M}(t)\rd_{H^1(\Omega)}\le C$ from Lemma~\ref{lem:KSM:arpriori} and 
 $0<\xi\le1$, the right-hand side of the preceding inequality is bounded by a constant independent of $\xi$. Therefore, taking the essential supremum over $t\in(0,T]$ yields the desired estimate.}
\end{proof}

\begin{proof}[Proof of Theorem~\ref{thm:conv:xi}]
{
Set $\widehat\xi_n=1/n$. Applying
Lemma~\ref{lem:nonlinear_compactness} to
$\{u^{\widehat\xi_n,M}\}_{n\in\N}$, we conclude that
$\{\beta(u^{\widehat\xi_n,M})\}_{n\in\N}$ is relatively compact in
$L^2(Q_T)$. Hence, Lemmas~\ref{lem:KSM:arpriori} and
\ref{lem:KSM:arpriori_v} yield a subsequence $\{\xi_k\}$ of
$\{\widehat\xi_n\}$ and functions
$u^*\in\mathcal{K}\cap H^1(0,T;H^1(\Omega)^*)$,
$\beta^*\in L^2(0,T;H^1(\Omega))$, and
$v_j^*\in L^2(0,T;H^1(\Omega))$ such that
}
{
\begin{align*}
&u^{\xi_k,M} \rightharpoonup u^*
& &  \mbox{weakly in } L^2(Q_T) \mbox{ and } H^1(0,T;H^1(\Omega)^*),  \\
&\beta \lr u^{\xi_k,M}\rr \to \beta^*
& & \mbox{strongly in } L^2(Q_T), \mbox{ a.e. in } Q_T,   \\
& & & \mbox{weakly in } L^2(0,T;H^1(\Omega)), \notag \\
&v_j^{\xi_k,M} \rightharpoonup v_j^* 
& & \mbox{weakly in } L^2(0,T;H^1(\Omega))  \notag 
\end{align*}
as $\xi_k$ tends to zero. 
}
{
Together with the weak convergence of $u^{\xi_k,M}$,
Lemma~6.1 in \citet{eymard1999fvm} gives
$\beta^*=\beta(u^*)$. Since $\beta$ is strictly increasing, its inverse
is continuous on $\beta([0,1])$, and hence
\[
u^{\xi_k,M}\to u^*
\quad\mbox{a.e. in }Q_T.
\]
The uniform bound $0\le u^{\xi_k,M}\le1$ then gives
\[
u^{\xi_k,M}\to u^*
\quad\mbox{strongly in }L^2(Q_T).
\]
}
In the weak form of the equation for $v_j^{\xi_k,M}$, i.e., 
\begin{align}
\xi_k \int_0^T \la \frac{\partial v_j^{\xi_k,M}}{\partial t}, \varphi \ra 
+d_j\int_0^T \la   \nabla v_j^{\xi_k,M}  ,\nabla \varphi \ra
 +  \int_0^T \la v_j^{\xi_k,M}, \varphi\ra
 = \int_0^T \la {u}^{\xi_k,M}, \varphi\ra \label{weak_vjxi}  
\end{align}
for $\varphi \in L^2(0,T;H^1(\Omega))$, { the time-derivative term vanishes as $k\to\infty$}. 
{
Indeed, Lemma~\ref{lem:KSM:arpriori_v} gives
\[
\xi_k\left\|\frac{\partial v_j^{\xi_k,M}}{\partial t}\right\|_{L^2(Q_T)}
\le C\xi_k^{1/4}\to0.
\]
Therefore, passing to the limit in $\xi_k$, we obtain
}
\begin{align}
d_j\int_0^T \la   \nabla v_j^*  ,\nabla \varphi \ra
 +  \int_0^T \la v_j^*, \varphi\ra
 = \int_0^T \la {u}^*, \varphi\ra  \label{weak_vjxi*}
\end{align}
together with $v_j^*\in L^2(0,T; H^2(\Omega))$. Hence, $v_j^*$ satisfies the equation for $v_j$ in \eqref{eq:KSPE}. 
Subtracting \eqref{weak_vjxi*} from \eqref{weak_vjxi} and choosing $\varphi = v_j^{\xi_k,M}-v_j^*$, we obtain 
\begin{equation*}
\begin{aligned}
   &d_j\ld \nabla v_j^{\xi_k,M}- \nabla v_j^*\rd_{L^2(Q_T)}^2
    + 
    \ld  v_j^{\xi_k,M}-v_j^* \rd_{L^2(Q_T)}^2\\
     &=
{-}\xi_k \int_0^T \la \frac{\partial v_j^{\xi_k,M}}{\partial t}, v_j^{\xi_k,M}-v_j^* \ra     
+
    \int_0^T \la {u}^{\xi_k,M} - {u}^*, v_j^{\xi_k,M}-v_j^* \ra \\
    &\le
\xi_k\ld \frac{\partial v_j^{\xi_k,M}}{\partial t} \rd_{L^2(Q_T)} \ld v_j^{\xi_k,M}-v_j^* \rd_{L^2(Q_T)}     
+
    \ld {u}^{\xi_k,M} - {u}^*\rd_{L^2(Q_T)} \ld v_j^{\xi_k,M}-v_j^* \rd_{L^2(Q_T)}. 
\end{aligned}
\end{equation*}
Therefore, thanks to Lemma~\ref{lem:KSM:arpriori_v}, we get the strong convergence of $v_j^{\xi_k,M}$ in $L^2(0,T;H^1(\Omega))$ as $\xi_k$ tends to zero. 
Furthermore, this leads to 
the convergence of $g\lr u^{\xi_k,M}\rr \nabla \sum_{j=1}^M a_j v_j^{\xi_k,M} 
\to 
g(u^*) \nabla \sum_{j=1}^M a_j v_j^*$ in $L^2(Q_T)$ as $\xi_k$ tends to zero.
Now, we can pass to the limit in $\xi_k$ in the  weak form of the equation for $u^{\xi_k,M}$. 
{
The weak convergence in $H^1(0,T;H^1(\Omega)^*)$, together with the
continuity of the trace operator, yields
\[
u^*(0)=u_0.
\]
}
Thereby,  we see that $\lr u^*, v_1^*, \dots, v_M^* \rr$ is the weak solution of \eqref{eq:KSPE}. 
\end{proof}


\section{Approximation of gradients of interaction potentials}
\label{sec:approx_kernels}

In this section, we show that { the gradient of a given radially symmetric interaction potential} can be {approximated by} a linear combination of the gradients of Green functions. In particular, to investigate the relationship between \eqref{eq:AD} and \eqref{eq:ADtilde}, we require an approximation { for the gradient of the interaction potential in $L^1(\Omega)$}. 
{
\citet{murakawa_tanaka2024} provided an approximation of integral kernels in the $C^1$ space on a one-dimensional periodic domain.
A key advantage of their result is that it provides a concrete method for determining the coefficients together with an error estimate.
More recently, \citet{IshiiTanaka2026} obtained constructive $W^{1,1}(\R^N)$-approximations for $N=1,2,3$ under stronger assumptions; see Remark~\ref{rem:constructive_approx}. 
In contrast, the approximation result established below is valid in arbitrary spatial dimensions under Assumption~\eqref{assum_W}, although it is nonconstructive and provides no error estimate.
}

{Here, we first show that the gradient of an arbitrary radial function can be approximated in $L^1(\R^N)$ by a linear combination of the gradients of the fundamental solutions, regardless of dimension.}
Define $k_j$ by 
\begin{eqnarray}
k_{j}(x) &:=& \dfrac{1}{d^{N/2}_j} G_{N}\left( \dfrac{|x|}{\sqrt{d_j}} \right), \label{def:kj} \\
G_{N}(|x|) &:=& \dfrac{1}{(2\pi)^{N/2}}\dfrac{1}{|x|^{N/2-1}} \mathcal{M}_{N/2-1}\left(|x|\right), 
\end{eqnarray}
where $\mathcal{M}_{\nu}$ is the modified Bessel function of the second kind {of order} $\nu$, namely, 
\[
\mathcal{M}_{\nu}(r) := \int^{+\infty}_{0} e^{-r \cosh s} \cosh(\nu s) ds. 
\]
{The asymptotic properties are} known { as follows; see, e.g., \citep{Watson}:}
\begin{equation}\label{mod:inf}
{{\mathcal{M}_{\nu}(r)}} \simeq\sqrt{\dfrac{\pi}{2r}} e^{-r}\quad (r\to+\infty)    
\end{equation}
and 
\begin{equation}\label{mod:zero}
    \mathcal{M}_{\nu}(r) \simeq 
    \begin{cases}
    \dfrac{\Gamma(\nu)}{2} \left(\dfrac{r}{2}\right)^{-\nu}, &(\nu>0),\\
    -\log r,  &(\nu=0),
    \end{cases}
    \quad ({r\to0^+})
\end{equation}
for any $\nu\ge 0$. 
{
It follows from \citep{Watson} that $\widehat{G_N}(\zeta)=(1+|\zeta|^2)^{-1}$. Hence, $G_N$ is the fundamental solution of $-\Delta+1$, and the scaling in \eqref{def:kj} shows that $k_j$ is the fundamental solution satisfying
}
\[
d_j\Delta k_j -k_j + \delta =0 \quad \mbox{ in } \R^N,
\]
where $\delta$ is the Dirac delta function.
{
Moreover, \eqref{mod:inf}--\eqref{mod:zero} imply that $k_j \in C\lr \R^N\setminus \{0\}\rr \cap W^{1,1}\lr \R^N \rr$.
}
Therefore, $v=k_j*u$ satisfies 
\[
d_j\Delta v -v + u =0 \quad \mbox{ in } \R^N. 
\]

We obtain the following result:
\begin{proposition}\label{thm:ker}
Let $K$ be a radial function satisfying $|\nabla K|\in L^1(\R^N)$.
Then, for any $\eps>0$, there exist $M$, a family of constants $\{a_j\}^M_{j=1}$, and a family of positive constants $\{d_j\}^M_{j=1}$ such that
\begin{equation*}
\left\| \nabla K - \nabla \sum^{M}_{j=1} a_j k_j \right\|_{L^1\lr \R^N\rr} < \eps
\end{equation*}
holds. 
\end{proposition}
\begin{proof}
    For a radial function $K$, we denote $K(x)=K_0(|x|)$.
    Then, the gradient is written as
    \begin{equation*}
        \nabla K(x) = \dfrac{x}{|x|} K'_0(|x|).
    \end{equation*}    
    Also, the gradient of the fundamental solution is written as
    \begin{equation*}
       \nabla k_j (x) = \dfrac{x}{d^{(N+1)/2}_j|x|} G'_{N}\left( \dfrac{|x|}{\sqrt{d_j}} \right). 
    \end{equation*}
    Then, for any given $M\in\N$, a family of constants $\{a_j\}^M_{j=1}$, and a family of positive constants $\{d_j\}^M_{j=1}$, we have
    \begin{eqnarray*}
    \left\|\nabla K - \nabla \sum^{M}_{j=1} a_j k_j \right\|_{L^1(\R^N)} 
    &=& |S^{N-1}| \int^{+\infty}_{0} r^{N-1} \left| K'_0(r) - \sum^{M}_{j=1} \dfrac{a_j}{d^{\frac{N+1}{2}}_j}  G'_{N}\left( \dfrac{r}{\sqrt{d_j}} \right) \right|dr \\
    &=& |S^{N-1}| \int_{\R} \left| e^{Ns} K'_0(e^{s}) - \sum^{M}_{j=1} b_j e^{N(s-h_j)} G'_{N}\left( e^{s-h_j} \right) \right|ds \\
    &=& |S^{N-1}| \int_{\R} \left| e^{Ns} K'_0(e^{s}) - \sum^{M}_{j=1} b_j F_N(s-h_j) \right|ds
    \end{eqnarray*}
    by replacing the variable $r=e^{s}$, $h_j = \frac{1}{2} \log d_j$, $b_j = a_j e^{-h_j}$,
    where $F_N(s):=e^{Ns} G'_N(e^s)$.
    {
    Moreover,
    \begin{equation*}
        \|F_N\|_{L^1(\R)}=\int_0^{+\infty}r^{N-1}|G'_N(r)|\,dr<+\infty.
    \end{equation*}
    Indeed, this follows from $G'_1(r)=-e^{-r}/2$ when $N=1$, while for $N\ge2$ it follows from
    \begin{equation*}
        G'_N(r)=-\dfrac{1}{(2\pi)^{N/2}}r^{1-N/2}\mathcal{M}_{N/2}(r)
    \end{equation*}
    and the asymptotic properties \eqref{mod:inf} and \eqref{mod:zero}.
    }
    {
    Here, the assumption $|\nabla K|\in L^1(\R^N)$ and the change of variables $r=e^s$ give
    \begin{equation*}
        \|\nabla K\|_{L^1(\R^N)}
        = |S^{N-1}| \int^{+\infty}_{0} r^{N-1}|K'_0(r)|\,dr
        = |S^{N-1}| \int_{\R} e^{Ns}|K'_0(e^{s})|\,ds
        <+\infty.
    \end{equation*}
    }
    Thus, it follows that by employing the variable transformations, the desired assertion is obtained by showing that $\mathrm{span}\{F_N(\cdot-h)\mid h\in\R\}$ is dense in $L^1(\R)$.

    According to Wiener's { Tauberian} theorem \citep{Wiener}, 
    $\mathrm{span}\{F_N(\cdot-h)\mid h\in\R\}$ is dense in $L^1(\R)$ if and only if the Fourier transform of $F_N$ has no real zero.
    When $N=1$, since $G'_1(r)=-e^{-r}/2$ holds, 
    we obtain
    {
    \begin{eqnarray*}
        \int_{\R} e^{i\zeta s} F_1(s) ds 
        &=& \int^{+\infty}_{0} \eta^{i\zeta} G'_1(\eta) d\eta 
        = -\dfrac{1}{2}\int^{+\infty}_{0} \eta^{i\zeta} e^{-\eta} d\eta 
        = - \dfrac{1}{2} \Gamma(1+i\zeta) \neq 0
    \end{eqnarray*}}
    { for any $\zeta\in\R$.}
    In the case $N\ge 2$, we use the integral formula \citep[\S~13.21]{Watson}
    \begin{equation*}
        \int^{+\infty}_{0} r^{\mu-1} \mathcal{M}_{\nu}(r)dr = 2^{\mu-2} \Gamma\left(\dfrac{\mu-\nu}{2}\right) \Gamma\left(\dfrac{\mu+\nu}{2}\right)
    \end{equation*}
    for $|\mathrm{Re}(\nu)|< \mathrm{Re}(\mu)$.
    {
    By setting $\mu=N/2+i\zeta$, $\nu=N/2-1$, it follows from substitution and integration by parts that 
    \begin{eqnarray*}
        \int_{\R} e^{i\zeta s} F_N(s)ds 
        &=& \int_{\R} e^{(N+i\zeta) s} G'_N(e^{s}) ds \\
        &=& \int^{+\infty}_{0} r^{N-1+i\zeta} G'_N(r) dr \\
        &=& -(N-1+i\zeta)\int^{+\infty}_{0} r^{N-2+i\zeta} G_N(r) dr \\
        &=& -\dfrac{N-1+i\zeta}{(2\pi)^{N/2}}\int^{+\infty}_{0} r^{N/2-1+i\zeta} \mathcal{M}_{N/2-1}(r) dr \\
        &=& -\dfrac{2^{-2+i\zeta}}{\pi^{N/2}} (N-1+i\zeta) \Gamma\left( \dfrac{1+i\zeta}{2} \right) \Gamma\left( \dfrac{N-1+i\zeta}{2} \right) \neq 0
    \end{eqnarray*}
    holds for any $\zeta\in\R$.
    }
    {
    Here, the boundary terms in the integration by parts vanish. Indeed, \eqref{mod:zero} gives
    \begin{equation*}
        r^{N-1}|G_N(r)|=
        \begin{cases}
            O(r|\log r|),&N=2,\\
            O(r),&N\ge3,
        \end{cases}
        \quad (r\to0),
    \end{equation*}
    while \eqref{mod:inf} gives
    \begin{equation*}
        r^{N-1}|G_N(r)|=O\left(r^{(N-1)/2}e^{-r}\right)
        \quad (r\to+\infty).
    \end{equation*}
    }
    Therefore, for any $N\in\N$, the Fourier transform of $F_N$ has no real zero.

    Since $\mathrm{span}\{F_N(\cdot-h)\mid h\in\R\}$ is dense in $L^1(\R)$, for any $\eps>0$,
    there exists $M\in\N$ and  families of constants $\{b_j\}^M_{j=1}$ and $\{h_j\}^M_{j=1}$ such that
    \beaa
    |S^{N-1}| \int_{\R} \left| e^{Ns} K'_0(e^{s}) - \sum^{M}_{j=1} b_j F_N(s-h_j) \right|ds < \eps
    \eeaa
    holds.
    Then, by setting $d_j = e^{2h_j},\ a_j=b_j e^{h_j}$ for { $j=1,\ldots,M$}, we obtain
    \begin{equation*}
        \left\|\nabla K - \nabla \sum^{M}_{j=1} a_j k_j \right\|_{L^1(\R^N)} <\eps
    \end{equation*}
    which completes the proof.
\end{proof}

{ We transfer the whole-space approximation in Proposition~\ref{thm:ker} to the periodic interaction field $\nabla W$. 
For a given radial function $K$ satisfying $|\nabla K|\in L^1(\R^N)$, we construct the periodic interaction field $\nabla W$ on $\Omega = [-L,L)^N = \R^N / (2L\Z)^N$ as in \eqref{wholespace_to_periodic}.
}
{
The series in \eqref{wholespace_to_periodic} converges absolutely in $L^1(\Omega)$ because
\begin{align}\label{peri2whole}
    \int_\Omega \left| \nabla W(x) \right| dx
    \le
     \sum_{\ell\in \Z^N}\int_\Omega \left| \nabla K(x-2L\ell) \right| dx
    =
    \int_{\R^N} \left| \nabla K(x) \right| dx .   
\end{align}
}

{
Since $k_j$ defined in \eqref{def:kj} belongs to $W^{1,1}\lr \R^N\rr$, the following series converges in $W^{1,1}(\Omega)$ and defines the Green function of the periodic problem.
}
\begin{equation*}
    {w_j:=\sum_{\ell\in \Z^N} k_j(\,\cdot-2L\ell)
    \quad \mbox{in }W^{1,1}(\Omega).}
\end{equation*}
We note that 
$v=w_j*u$ solves  
\[
d_j\Delta v -v + u =0 \quad \mbox{ in } \Omega. 
\]

\begin{proof}[Proof of Theorem~\ref{thm:W_approx}]
Assume that $W$ satisfies \eqref{assum_W}. Then, \eqref{wholespace_to_periodic} holds for a certain radial function $K$ satisfying $|\nabla K|\in L^1(\R^N)$. 
By Proposition~\ref{thm:ker}, for any $\eps > 0$, there exist some $M \in \N$, $\{a_j\}^M_{j=1}$ and $\{d_j\}^M_{j=1}$ such that 
\begin{equation*}
\left\| \nabla K - \nabla \sum^{M}_{j=1} a_j k_j \right\|_{L^{1}\lr \R^N\rr} < \eps.
\end{equation*}
Then, it follows from \eqref{peri2whole} that 
\begin{align}
    \left\| \nabla W - \nabla \sum^{M}_{j=1} a_j w_j \right\|_{L^{1}\lr \Omega\rr}
    \le
    \left\| \nabla K - \nabla \sum^{M}_{j=1} a_j k_j \right\|_{L^{1}\lr \R^N \rr}
    <\eps. \label{WM-barWM}
\end{align}
holds. The proof is complete. 
\end{proof}

{
\begin{remark}\label{rem:constructive_approx} 
Under stronger regularity and decay assumptions, \citet{IshiiTanaka2026} established constructive approximations of radial kernels by finite linear combinations of fundamental solutions in $W^{1,1}(\R^N)$ for $N=1,2,3$.
More precisely, they fixed the diffusion coefficients as $d_j=j^{-2}$, determined the coefficients $a_j$ through polynomial approximations, and derived quantitative error estimates; see \cite[Section~4]{IshiiTanaka2026}. 

Using \eqref{peri2whole}, their whole-space result provides a constructive counterpart to Theorem~\ref{thm:W_approx} for $N=1,2,3$ under the additional assumptions imposed in \cite{IshiiTanaka2026}. 
Thus, explicit coefficients and quantitative error estimates are available for the corresponding restricted class of periodic interaction fields. 
This does not make Theorem~\ref{thm:W_approx} constructive in its full generality. 
\end{remark}
}


{
\section{Convergence under kernel approximation and the main comparison result}
\label{sec:convWm}
}

{ In this section, we first prove Theorem~\ref{thm:conv:Wm}.
We then deduce Theorem~\ref{thm:exist:AD} and finally establish the main result, Theorem~\ref{Thm_conv}.}

\begin{proof}[Proof of Theorem~\ref{thm:conv:Wm}]
{
Since $\nabla W_m\to\nabla W$ in $L^1(\Omega)$,
$\{\nabla W_m\}_{m\in\N}$ is bounded in $L^1(\Omega)$.
Set
\[
\Phi(s):=\int_0^s\beta(r)\,dr.
\]
Since
\[
\beta(u_m)\in L^2(0,T;H^1(\Omega)),
\qquad
\frac{\partial u_m}{\partial t}\in L^2(0,T;H^1(\Omega)^*),
\]
the nonlinear chain rule of
\citet[Lemma~1.5]{altLuckhaus1983}, applied with
$w_m=\beta(u_m)$ and $b=\beta^{-1}$ on $\beta([0,1])$, yields
\[
\int_0^t
\left\langle\frac{\partial u_m}{\partial t},\beta(u_m)\right\rangle ds
=
\int_\Omega\Phi(u_m(t))\,dx
-
\int_\Omega\Phi(u_0)\,dx
\]
for a.e.\ $t\in(0,T)$.
Testing the equation for $u_m$ by $\beta(u_m)$ and estimating
its time derivative in $H^1(\Omega)^*$ as in
Section~\ref{sec:existence}, we obtain
\[
\sup_{m\in\N}\left(
\left\|u_m\right\|_{H^1(0,T;H^1(\Omega)^*)}
+
\left\|\beta(u_m)\right\|_{L^2(0,T;H^1(\Omega))}
\right)<+\infty.
\]
Applying Lemma~\ref{lem:nonlinear_compactness} to $\{u_m\}_{m\in\N}$,
we obtain a subsequence $\{m_k\}_{k\in\N}$,
$u_*\in\mathcal{K}\cap H^1(0,T;H^1(\Omega)^*)$, and
$\beta^*\in L^2(0,T;H^1(\Omega))$ such that
\begin{align*}
&u_{m_k} \rightharpoonup u_*
& &  \mbox{weakly in } L^2(Q_T)
\mbox{ and } H^1(0,T;H^1(\Omega)^*), \\
&\beta \lr u_{m_k}\rr \to \beta^*
& & \mbox{strongly in } L^2(Q_T), \mbox{ a.e. in } Q_T, \\
&\beta \lr u_{m_k}\rr \rightharpoonup \beta^*
& & \mbox{weakly in } L^2(0,T;H^1(\Omega))
\end{align*}
as ${m_k}$ tends to infinity. 
Together with the weak convergence of $u_{m_k}$,
Lemma~6.1 in \citet{eymard1999fvm} gives
$\beta^*=\beta(u_*)$. Since $\beta$ is strictly increasing, its
inverse is continuous on $\beta([0,1])$, and hence
$u_{m_k}\to u_*$ a.e. in $Q_T$. The uniform bound
$0\le u_{m_k}\le1$ then implies
\[
u_{m_k}\to u_*
\quad\mbox{strongly in }L^2(Q_T).
\]
}
To take the limit in the weak form of the equation, evaluating the convergence of $\nabla W_{m_k} * u_{m_k}$ is necessary. 
However, this can be easily obtained using Young's inequality as follows: 
\begin{align}
    &\ld \nabla W_{m_k} * u_{m_k} -\nabla W * u_{*} \rd_{L^{2}\lr Q_T \rr}\notag\\
    &\quad \le
    \ld \lr \nabla W_{m_k}-\nabla W\rr * u_{m_k} \rd_{L^{2}\lr Q_T \rr}
    +
    \ld \nabla W * \lr u_{m_k} - u_{*}\rr \rd_{L^{2}\lr Q_T \rr}\notag\\
    &\quad \le
    { \ld  \nabla W_{m_k}-\nabla W \rd_{L^{1}(\Omega)}} \ld u_{m_k} \rd_{L^{2}\lr Q_T \rr}
    +
    { \ld \nabla W \rd_{L^{1}(\Omega)}} \ld  u_{m_k} - u_{*} \rd_{L^{2}\lr Q_T \rr}\notag\\
    & \qquad \to 0 \mbox{ as } m_k \to \infty. \label{eq:drift-field-convergence}
\end{align}
We also have the strong convergence of $g(u_{m_k}) \nabla W_{m_k} * u_{m_k}$ to $g(u_{*}) \nabla W * u_{*}$ in ${L^{2}\lr Q_T \rr}$. 
{
The weak convergence in $H^1(0,T;H^1(\Omega)^*)$, together with the
continuity of the trace operator, yields
\[
u_*(0)=u_0.
\]
}
Therefore, passing to the limit in $m_k$ in the weak form of \eqref{eq:ADtilde}, the desired result holds. 
\end{proof}

{ Next, we establish the existence of a weak solution to \eqref{eq:AD}.}

\begin{proof}[Proof of Theorem~\ref{thm:exist:AD}]
    It follows from Theorem~\ref{thm:W_approx} that there exists a sequence of functions $\{ W_m\}_{m\in \N}$ such that $\nabla W_m$ converges to $\nabla W$ in $L^{1}(\Omega)$ as $m$ tends to infinity.
    { More precisely, for each $m$, we write
    \[
    W_m=\sum_{j=1}^{M_m}a_j^{(m)}w_j^{(m)},
    \]
    where $w_j^{(m)}$ is the Green function corresponding to $d_j^{(m)}$. Since $w_j^{(m)}\in W^{1,1}(\Omega)$, Young's inequality gives
    \[
    \|w_j^{(m)}*u_0\|_{W^{1,\infty}(\Omega)}
    \le
    \|w_j^{(m)}\|_{W^{1,1}(\Omega)}
    \|u_0\|_{L^\infty(\Omega)}.
    \]
    Thus the compatible initial data satisfy Assumption~\eqref{assum_init_v}.}
    { Hence, Theorems~\ref{thm:exist:KSM} and \ref{thm:conv:xi}, applied to \eqref{eq:KSM} with the compatible initial data \eqref{init_u_wu}, yield a weak solution of \eqref{eq:KSPE} satisfying $u(0)=u_0$, which, in turn, means the existence of a weak solution to \eqref{eq:ADtilde}.} 
    Thus, Theorem~\ref{thm:conv:Wm} establishes the existence of a weak solution of \eqref{eq:AD}. 
\end{proof}

{
Finally, we prove the main comparison result, Theorem~\ref{Thm_conv}.

\begin{proof}[Proof of Theorem~\ref{Thm_conv}]
Let $\eps>0$ and $T>0$ be fixed. By Theorem~\ref{thm:W_approx}, for each $n\in\N$, we can choose
$
M_n\in\N,
\{a_{j}^{(n)}\}_{j=1}^{M_n},
\{d_{j}^{(n)}\}_{j=1}^{M_n}$,
let $w_{j}^{(n)}$ be the Green function corresponding to $d_{j}^{(n)}$, and set
\[
W_n=\sum_{j=1}^{M_n}a_{j}^{(n)}w_{j}^{(n)}.
\]
These quantities are chosen so that
\[
\ld\nabla W_n-\nabla W\rd_{L^1(\Omega)}<\frac{1}{n}.
\]

For each $n$, choose a positive sequence $\{\xi_\ell^{(n)}\}_{\ell\in\N}$ such that $\xi_\ell^{(n)}\to0$ and corresponding weak solutions of the parabolic--parabolic system~\eqref{eq:KSM} with the compatible initial data~\eqref{init_u_wu}. By Theorems~\ref{thm:exist:KSM} and \ref{thm:conv:xi}, along a subsequence of $\ell$, these weak solutions converge to a weak solution
\[
\lr u_n,v_{n,1},\ldots,v_{n,M_n}\rr
\]
of the parabolic--elliptic system~\eqref{eq:KSPE}. In particular, $u_n$ is a weak solution of \eqref{eq:ADtilde} corresponding to the kernel $W_n$.

Applying Theorem~\ref{thm:conv:Wm} to the sequence $\{u_n\}_{n\in\N}$ and taking a subsequence if necessary, we obtain a weak solution $u$ of \eqref{eq:AD} such that
\[
u_{n_k}\to u\qquad\mbox{strongly in }L^2(Q_T).
\]
We fix a sufficiently large $n_k$ such that
\[
\ld u-u_{n_k}\rd_{L^2(Q_T)}<\frac{\eps}{2}.
\]
Set
\[
M:=M_{n_k},
\]
and define
\[
\{a_j\}_{j=1}^{M}:=
\{a_{j}^{(n_k)}\}_{j=1}^{M_{n_k}},
\qquad
\{d_j\}_{j=1}^{M}:=
\{d_{j}^{(n_k)}\}_{j=1}^{M_{n_k}}.
\]
Since $u_{n_k}$ was selected as a subsequential limit in Theorem~\ref{thm:conv:xi}, we can choose a sufficiently large $\ell$ along the corresponding subsequence so that
\[
\ld u_{n_k}-u^{\xi_\ell^{(n_k)},M}\rd_{L^2(Q_T)}<\frac{\eps}{2}.
\]
Setting $\xi:=\xi_\ell^{(n_k)}$ and applying the triangle inequality, we obtain
\[
\ld u-u^{\xi,M}\rd_{L^2(Q_T)}<\eps.
\]
Thus, we have the desired result. \end{proof}
}

\section{Numerical simulations}\label{sec:numerical}
{\red
In this section, we numerically investigate the relationship between the nonlocal haptotaxis model and the Keller--Segel type chemotaxis systems on a one-dimensional periodic domain. We first examine finite-sum approximations of an interaction potential and its gradient. We then compare numerical solutions of the nonlocal model, the corresponding parabolic--elliptic system, and the parabolic--parabolic system in order to investigate how the number $M$ of terms in the kernel approximation and the relaxation time $\xi$ affect the differences among the numerical solutions.


{
The results of this paper establish the approximation of the gradient of an interaction potential by finite linear combinations of gradients of Green functions in arbitrary spatial dimensions. 
Under the general assumption~\eqref{assum_W}, however, the proof based on Wiener's Tauberian theorem does not provide an explicit construction of the coefficients and diffusion constants. 
Constructive approximations for $N=1,2,3$ under stronger assumptions are available from \citet{IshiiTanaka2026}, as discussed in Remark~\ref{rem:constructive_approx}. 

The numerical study below is performed on a one-dimensional periodic
domain in order to isolate the effects of the number $M$ of kernel terms and the relaxation time $\xi$ in a simple
and reproducible setting.
We use the construction of \citet{murakawa_tanaka2024} because it is formulated directly on a one-dimensional periodic domain and provides the periodic Green functions and coefficients used below without an additional periodization step.
We note that, in one spatial dimension, the construction of \citet{IshiiTanaka2026} is closely related to the periodic construction of \citet{murakawa_tanaka2024}: both are based on the choice $d_j=j^{-2}$, apart from differences in the indexing and the treatment of the constant mode, and on polynomial approximation using Chebyshev nodes.
}

Let $\Omega=[-4,4)$ and impose the periodic boundary condition. As an interaction potential, we use the following continuously differentiable, compactly supported raised-cosine kernel:
\begin{equation}\label{eq:numerical-kernel}
W(x)=
\begin{cases}
\displaystyle \frac{1}{2}\bigl(1+\cos(\pi x)\bigr), & |x|<1,\\[1mm]
0, & |x|\geq 1.
\end{cases}
\end{equation}
Following \citet{murakawa_tanaka2024}, we set $d_j=1/j^2$ and approximate $W$ by
\[
W_M(x)=a_0^{(M)}+\sum_{j=1}^M a_j^{(M)}w_j(x),
\qquad
w_j(x)=
\frac{\cosh\bigl((4-|x|)/\sqrt{d_j}\bigr)}
{2\sqrt{d_j}\sinh\bigl(4/\sqrt{d_j}\bigr)}.
\]
The coefficients $\{a_j^{(M)}\}_{j=0}^M$ are computed separately for each $M$ by degree-$M$ Lagrange interpolation at $M+1$ Chebyshev--Gauss nodes. More precisely, after the change of variable $y=\cosh(4-|x|)$, which maps $0\leq |x|\leq4$ onto $1\leq y\leq\cosh4$, we use
\[
y_\ell=\frac{1+\cosh4}{2}
+\frac{\cosh4-1}{2}
\cos\left(\frac{(2\ell+1)\pi}{2(M+1)}\right),
\qquad \ell=0,\ldots,M.
\]
Thus, the interpolation values are $W(4-\operatorname{arcosh}y_\ell)$. Using the identity $T_j(\cosh z)=\cosh(jz)$ for the Chebyshev polynomial $T_j$, the resulting degree-$M$ polynomial in $y$ is then converted into a linear combination of the Green functions $w_j$. To measure the approximation of the interaction kernel, we define the relative error
\[
e_W(M)=
\frac{\|\partial_xW-\partial_xW_M\|_{L^1(\Omega)}}
{\|\partial_xW\|_{L^1(\Omega)}}.
\]

Figure~\ref{fig:kernel-approximation} compares $W$ and $\partial_xW$ with their approximations for $M=3$ and $M=13$, and shows the dependence of $e_W(M)$ on $M$.
\begin{figure}[bt]
\begin{center}
\includegraphics[width=\textwidth]{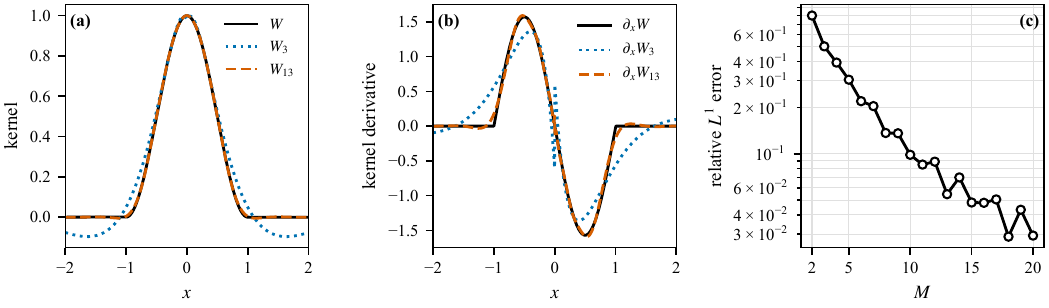}
\end{center}
\caption{{Approximation of the raised-cosine interaction potential~\eqref{eq:numerical-kernel}. (a) Profiles of $W$, $W_3$, and $W_{13}$. (b) Profiles of $\partial_xW$, $\partial_xW_3$, and $\partial_xW_{13}$. (c) Relative $L^1$ error $e_W(M)$ of the kernel derivatives for $M=2,\ldots,20$}}
\label{fig:kernel-approximation}
\end{figure}

For the time-evolution simulations, on the interval $0\leq u\leq1$, we set
\[
\chi(u)=u^3,\qquad
\beta(u)=\frac{3}{4}u^4,\qquad
g(u)=u(1-u^3).
\]
Outside this interval, $\beta$ and $g$ are understood to be extended consistently with Assumptions~\eqref{assum_gamma} and \eqref{assum_g}; in particular, $g(u)=0$ for $u\notin[0,1]$. The computed solutions remain in $[0,1]$, so these extensions do not affect the simulations.
The final time, the number of spatial cells, and the spatial and temporal mesh sizes are
\[
T=100,\qquad N_x=300,\qquad
\Delta x=\frac{8}{300},\qquad
\Delta t=5\times10^{-5},
\]
respectively. The periodic cell centers are denoted by $x_i=-4+i\Delta x$, $i=0,\ldots,N_x-1$. To compare the aggregation dynamics emerging from the homogeneous state $u=1/2$, we impose a deterministic small perturbation with modal amplitude $0.01$. The first five Fourier modes are included so that the comparison does not depend on a single spatial frequency. We use cosine functions for the odd modes and sine functions for the even modes to avoid an evolution constrained by reflection symmetry. This choice allows the three models to be compared under the same reproducible initial condition without dependence on a random perturbation. We use the initial values
\[
\begin{aligned}
u_i^0=\frac12+0.01\bigg\{&
\cos\left(\frac{\pi x_i}{4}\right)
+\sin\left(\frac{\pi x_i}{2}\right)
+\cos\left(\frac{3\pi x_i}{4}\right)\\
&+\sin(\pi x_i)
+\cos\left(\frac{5\pi x_i}{4}\right)
\bigg\}.
\end{aligned}
\]
For the parabolic--parabolic system, the compatible initial values corresponding to \eqref{init_u_wu} are given explicitly by
\[
\begin{aligned}
v_{j,i}^0=\frac12+0.01\bigg\{&
\frac{\cos(\pi x_i/4)}{1+d_j\lambda_{1,h}}
+\frac{\sin(\pi x_i/2)}{1+d_j\lambda_{2,h}}
+\frac{\cos(3\pi x_i/4)}{1+d_j\lambda_{3,h}}\\
&+\frac{\sin(\pi x_i)}{1+d_j\lambda_{4,h}}
+\frac{\cos(5\pi x_i/4)}{1+d_j\lambda_{5,h}}
\bigg\},
\end{aligned}
\]
where
\[
\lambda_{r,h}
=\frac{4}{\Delta x^2}
\sin^2\left(\frac{r\pi\Delta x}{8}\right),
\qquad r=1,\ldots,5,
\]
is the eigenvalue of the negative discrete Laplacian associated with the $r$th Fourier mode.

For the nonlocal model, the convolution at each time step is computed by the periodic quadrature
\[
\Psi_i^n
=\Delta x\sum_{\ell=0}^{N_x-1}
W(x_i-x_\ell)u_\ell^n,
\]
where $x_i-x_\ell$ is interpreted as the periodic distance. To describe the common population update for the three models, set
\[
\phi_i^n=
\begin{cases}
\Psi_i^n,&\text{for the nonlocal model},\\
\displaystyle\sum_{j=1}^M a_j^{(M)}v_{j,i}^n,
&\text{for the parabolic--elliptic and parabolic--parabolic systems}.
\end{cases}
\]
The local velocity at the interface $x_{i+1/2}$ is
\[
c_{i+1/2}^n=\frac{\phi_{i+1}^n-\phi_i^n}{\Delta x}.
\]
The population equations in the three models are discretized in space by the same finite-volume method. The numerical flux is
\[
\begin{aligned}
J_{i+1/2}^n={}&-\frac{\beta(u_{i+1}^n)-\beta(u_i^n)}{\Delta x}
+\frac{c_{i+1/2}^n}{2}\bigl(g(u_i^n)+g(u_{i+1}^n)\bigr)\\
&-\frac{\alpha_{i+1/2}^n}{2}(u_{i+1}^n-u_i^n),
\end{aligned}
\]
where the local Lax--Friedrichs speed is
\[
\alpha_{i+1/2}^n
=|c_{i+1/2}^n|
\max\bigl\{|g'(u_i^n)|,|g'(u_{i+1}^n)|\bigr\}.
\]
The population density is advanced by the forward Euler update
\[
u_i^{n+1}=u_i^n-\frac{\Delta t}{\Delta x}
\bigl(J_{i+1/2}^n-J_{i-1/2}^n\bigr).
\]
The time step $\Delta t$ is chosen to satisfy the stability restriction for this explicit update. The elliptic equations in the parabolic--elliptic system are discretized by centered finite differences. For the chemical equations in the parabolic--parabolic system, we use centered finite differences in space and a first-order semi-implicit Euler method in time, treating the diffusion and degradation terms implicitly and the source term explicitly.

Figure~\ref{fig:solution-profiles} compares the numerical solution $u$ of the nonlocal model with the numerical solutions $u^{\xi,M}$ of the parabolic--parabolic systems corresponding to $(\xi,M)=(2^{-9},3)$, $(2^{-1},13)$, and $(2^{-9},13)$ at $t=0$, $10$, $20$, $40$, $70$, and $100$. For $M=3$, the solution differs substantially from that of the nonlocal model and develops multiple peaks even when $\xi=2^{-9}$ is small. As observed in Figure~\ref{fig:kernel-approximation}(b), $\partial_xW_3$ contains an additional repulsive component at intermediate distances that is absent from the target interaction kernel. The formation of multiple peaks is consistent with this discrepancy in the kernel approximation.

In contrast, for $M=13$ and $\xi=2^{-9}$, the parabolic--parabolic solution closely follows the nonlocal solution throughout the computation. When $M=13$ and $\xi=2^{-1}$, a difference between the two solutions is visible at intermediate times. This difference reflects the delayed response caused by the finite relaxation time of the chemical fields. Here, $M$ parametrizes the spatial approximation of the interaction kernel, whereas $\xi$ parametrizes the response time of the chemical fields.

\begin{figure}[bt]
\begin{center}
\includegraphics[width=\textwidth]{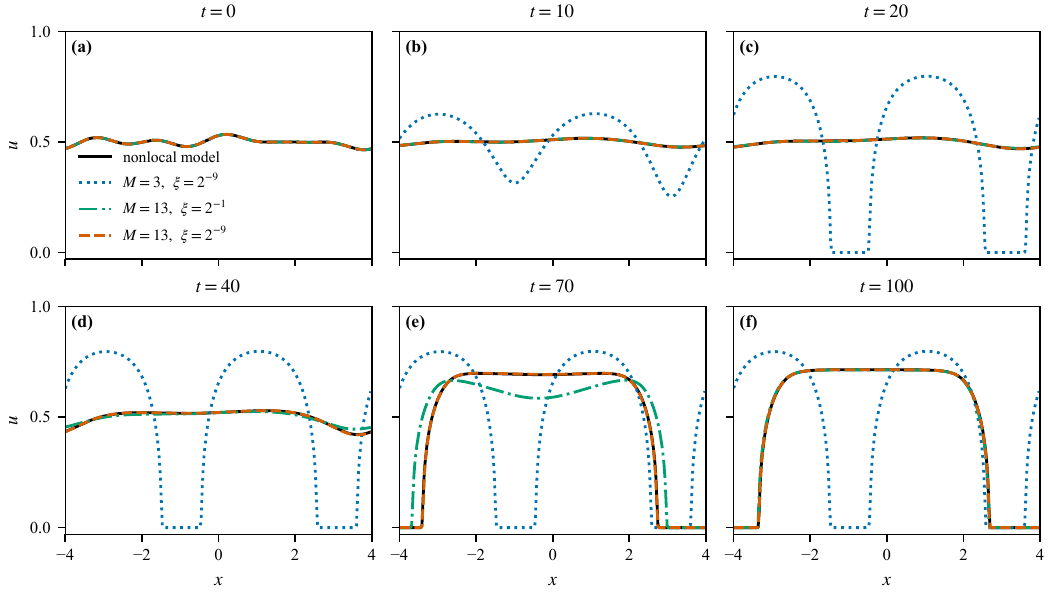}
\end{center}
\caption{{Numerical solutions of the nonlocal model and the parabolic--parabolic systems at $t=0$, $10$, $20$, $40$, $70$, and $100$. The black solid curves represent $u$, while the blue dotted, green dash-dotted, and orange dashed curves represent $u^{\xi,M}$ for $(\xi,M)=(2^{-9},3)$, $(2^{-1},13)$, and $(2^{-9},13)$, respectively}}
\label{fig:solution-profiles}
\end{figure}

Let $u$, $u^M$, and $u^{\xi,M}$ denote the numerical solutions of the nonlocal model, the parabolic--elliptic system, and the parabolic--parabolic system, respectively. All three models are computed on the same spatial grid with the same time step. For evaluating the relative differences, the numerical solutions are stored at $t_n=n\tau$, where $\tau=0.1$, $n=0,\ldots,N_o$, and $N_o=T/\tau$. We use the discrete space--time norm
\[
\|z\|_{L_h^2(Q_T)}^2
=\Delta x\,\tau
\left\{
\frac12\sum_{i=0}^{N_x-1}|z_i^0|^2
+\sum_{n=1}^{N_o-1}\sum_{i=0}^{N_x-1}|z_i^n|^2
+\frac12\sum_{i=0}^{N_x-1}|z_i^{N_o}|^2
\right\},
\]
which applies the rectangle rule in space and the trapezoidal rule in time. We define the following three relative differences between the numerical solutions:
\[
E_1(M)
=\frac{\|u-u^M\|_{L_h^2(Q_T)}}
{\|u\|_{L_h^2(Q_T)}},
\]
\[
E_2(M,\xi)
=\frac{\|u^M-u^{\xi,M}\|_{L_h^2(Q_T)}}
{\|u\|_{L_h^2(Q_T)}},
\]
\[
E_3(M,\xi)
=\frac{\|u-u^{\xi,M}\|_{L_h^2(Q_T)}}
{\|u\|_{L_h^2(Q_T)}}.
\]
Here, $E_1$ measures the difference between the nonlocal and parabolic--elliptic solutions, $E_2$ measures the difference between the parabolic--elliptic and parabolic--parabolic solutions, and $E_3$ measures the difference between the nonlocal and parabolic--parabolic solutions. The use of the common denominator gives
\[
E_3(M,\xi)\leq E_1(M)+E_2(M,\xi).
\]

Figure~\ref{fig:solution-errors}(a) shows $E_1(M)$ for $M=3$, $5$, $7$, $9$, $11$, and $13$. Figure~\ref{fig:solution-errors}(b) shows $E_2(M,\xi)$ for $M=5$, $9$, and $13$ with
\[
\xi=2^{-k},\qquad k=1,\ldots,9.
\]
Figure~\ref{fig:solution-errors}(c) shows $E_3(M,\xi)$ for $\xi=2^{-1}$, $2^{-5}$, $2^{-7}$, and $2^{-9}$, together with $E_1(M)$.

In the present numerical example, $E_1(M)$ decreases over the selected values of $M$, and $E_2(M,\xi)$ decreases as $\xi$ is reduced for each of the three values of $M$. Moreover, when $\xi$ is sufficiently small, $E_3(M,\xi)$ approaches $E_1(M)$, showing that the difference between the nonlocal and parabolic--parabolic solutions is then close to that between the nonlocal and parabolic--elliptic solutions.

We emphasize, however, that the coefficient sets $\{a_j^{(M)}\}_{j=0}^M$ are computed independently for each $M$ and do not form a nested sequence of approximations. Consequently, neither the kernel approximation error nor the difference between the numerical solutions is guaranteed to decrease monotonically with $M$. Local variations are in fact visible in Figures~\ref{fig:kernel-approximation}(c) and \ref{fig:solution-errors}(c). The plots describe only the behavior observed in the present numerical example and should not be interpreted as establishing a general quantitative relationship with respect to either $M$ or $\xi$.

\clearpage
\begin{figure}[bt]
\begin{center}
\includegraphics[width=\textwidth]{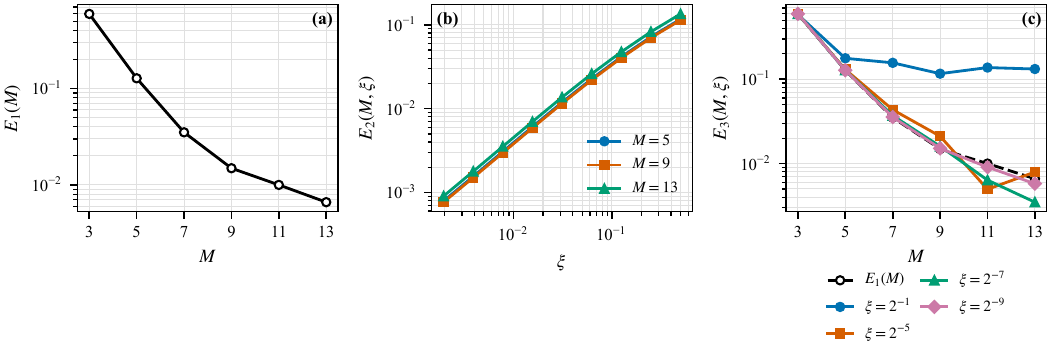}
\end{center}
\caption{{Relative differences between the numerical solutions. (a) $E_1(M)$ between the nonlocal and parabolic--elliptic solutions. (b) $E_2(M,\xi)$ between the parabolic--elliptic and parabolic--parabolic solutions for $M=5$, $9$, and $13$. (c) $E_3(M,\xi)$ between the nonlocal and parabolic--parabolic solutions for $\xi=2^{-1}$, $2^{-5}$, $2^{-7}$, and $2^{-9}$, together with $E_1(M)$}}
\label{fig:solution-errors}
\end{figure}
}

\section{Conclusion}
{
In this study, we investigated the relationship between haptotaxis and chemotaxis at the level of PDE models for cell sorting. In the nonlocal haptotaxis model, $\nabla W(z)$ is the interaction kernel that encodes the direction and magnitude of the contribution from cells at relative position $z$. Accordingly, $\nabla W*u$ represents the net interaction field obtained by integrating the surrounding cell density with respect to distance and direction. Thus, information acquired through cell contact or cellular protrusions is represented directly in terms of the cell density $u$, without introducing an intermediate field. By contrast, in the chemotaxis model, chemical substances produced by cells diffuse and degrade to form the concentration fields $v_j$, and the cells respond to their gradients $\nabla v_j$. Hence, cell--cell interactions are described directly by the kernel $\nabla W$ in the haptotaxis model, whereas they are transmitted indirectly through chemical fields in the chemotaxis model.

In the quasi-static chemotaxis system~\eqref{eq:KSPE}, each chemical concentration satisfies $v_j=w_j*u$. Consequently, the kernel approximation~\eqref{WM-barWM}, together with the estimate~\eqref{eq:drift-field-convergence} based on Young's inequality, shows that the haptotactic interaction field $\nabla W*u$ can be approximated by the superposition of chemical concentration gradients $\nabla\sum_{j=1}^M a_jv_j$. Here, the coefficients $a_j$ and $d_j$ are parameters in the kernel approximation and should not be interpreted as uniquely identified chemical substances or signalling pathways. The temporal scale is controlled separately by the common relaxation parameter $\xi$. We do not claim that the present results describe specific biological processes such as receptor kinetics or cytoskeletal dynamics.

We also numerically examined the differences among the numerical solutions of the three models on a one-dimensional periodic domain when the number $M$ of terms in the kernel approximation and the relaxation time $\xi$ were varied. For $M=3$, the approximated interaction kernel contains a repulsive component at intermediate distances, and the numerical solution developed multiple peaks that were not observed in the nonlocal model. In contrast, for $M=13$ and sufficiently small $\xi$, the difference between the numerical solutions of the nonlocal model and the parabolic--parabolic chemotaxis system was small.

In the present numerical example, the relative difference $E_1(M)$ between the numerical solutions of the nonlocal model and the parabolic--elliptic system decreased over the selected values of $M$. The relative difference $E_2(M,\xi)$ between the parabolic--elliptic and parabolic--parabolic systems also decreased as $\xi$ was reduced. Furthermore, for sufficiently small $\xi$, the relative difference $E_3(M,\xi)$ between the nonlocal and parabolic--parabolic models was close to $E_1(M)$. However, the coefficient sets $\{a_j^{(M)}\}$ are computed independently for each $M$ and do not form a nested sequence. Therefore, neither the kernel approximation error nor the difference between the numerical solutions is guaranteed to decrease monotonically with $M$. These observations are specific to the present numerical example and do not establish a general quantitative relationship with respect to $M$ or $\xi$.

Finally, our analysis also provides a mathematical contribution: we show that the gradient of a given radially symmetric interaction potential can be approximated by linear combinations of gradients of fundamental solutions in arbitrary spatial dimensions, and that suitable weak solutions of the corresponding models can be chosen so that their $L^2(Q_T)$ difference is arbitrarily small over a finite time interval. Since these are finite-time results, they do not address long-time pattern selection, for which the two models may differ. Further progress will require quantitative error estimates, long-time analysis, and extensions to multiple cell types and their interactions.
}


\section*{Acknowledgments}
{ The authors would like to thank the anonymous reviewers for their careful reading of the manuscript and their valuable comments and suggestions.}
The authors were partially supported by JSPS KAKENHI Grant Numbers 24H00188 and 22K03444.
HM was partially supported by JSPS KAKENHI Grant Number 21KK0044. 
YT was partially supported by JSPS KAKENHI Grant Number 24K06848.
{HI was partially supported by JSPS KAKENHI Grant Numbers 23K13013 and 25K07142.}


{
\section*{Data availability}
The numerical results reported in this article were generated using
computer code developed by the authors. The code is available from the
corresponding author upon reasonable request. No external datasets were
used in this study.
}

\section*{Conflict of interest}
The authors declare that they have no conflict of interest. 


\end{document}